\documentclass[11pt,leqno]{article}

\usepackage{graphicx}
\usepackage{latexsym} 
\usepackage{amsmath}
\usepackage{amssymb}
\usepackage{amsfonts}
\usepackage{enumerate}
\usepackage{theorem}

\theoremstyle{plain}
\newtheorem{teo}{Theorem}[section]
\newtheorem{lem}[teo]{Lemma}
\newtheorem{cor}[teo]{Corollary}
\newtheorem{prop}[teo]{Proposition}

{\theorembodyfont{\rmfamily}
\newtheorem{defin}[teo]{Definition}
\newtheorem{oss}[teo]{Remark}

\newtheorem{assum}[teo]{Assumption}
}

\renewcommand{\eqref}[1]{\textnormal{(\ref{#1})}}

\numberwithin{equation}{section}

\newcommand{\cvd}{\hfill$\square$}
\newcommand{\proof}[1]{\noindent\textsc{Proof#1}}

\title{Interior decay of solutions to elliptic equations with respect to frequencies at the boundary}

\author{Michele Di Cristo\thanks{Dipartimento di Matematica, Politecnico di Milano, Italy.\newline E-mail:
\texttt{michele.dicristo@polimi.it}} \and\
Luca Rondi\thanks{Dipartimento di Matematica,
Universit\`a degli Studi di Milano, Italy.\newline E-mail: \texttt{luca.rondi@unimi.it}} }

\date{}

\begin{document}

\maketitle

\setcounter{section}{0}
\setcounter{secnumdepth}{2}

\begin{abstract}
We prove decay estimates in the interior
for solutions to elliptic equations in divergence form with Lipschitz continuous coefficients. The estimates explicitly depend on the distance from the boundary and on suitable notions of frequency of the Dirichlet boundary datum. We show that, as the frequency at the boundary grows, the square of a suitable norm of the solution in a compact subset of the domain decays in an inversely proportional manner with respect to the corresponding frequency.

Under Lipschitz regularity assumptions, these estimates are essentially optimal and they have important consequences for the choice of optimal measurements for corresponding inverse boundary value problem.

\smallskip



\smallskip

\end{abstract}

\section{Introduction}\label{intro}

An important motivation for our study comes from elliptic inverse boundary value problems, such as the Calder\'on problem. Let us consider a bounded domain $\Omega\subset \mathbb{R}^N$, $N\geq 2$, which is regular enough and let $\gamma$ be a positive bounded function which is bounded away from zero and that corresponds to the background conductivity of a conducting body contained in $\Omega$. The aim of the inverse problem is to recover perturbations of the background conductivity, for example inhomogeneities, by performing suitable electrostatic measurements at the boundary of current and voltage type.
Such a problem comes from several types of nondesctructive evaluation problems in materials, where the aim is to detect the presence of flaws, as well as from medical imaging problems, where the aim is to detect the presence of tumors.

Namely, if $\tilde{\gamma}$ is the perturbed conductivity, one usually prescribes the voltage $f$ on the boundary of $\Omega$ and measures the corresponding current still on the boundary, that is,
$\tilde{\gamma}\nabla \tilde{u}\cdot \nu$ on $\partial\Omega$, where $\nu$ is the outer normal on $\partial\Omega$ and $\tilde{u}$, the electrostatic potential in $\Omega$, is the solution to the Dirichlet boundary value problem
\begin{equation}\label{dirpbmpert}
\left\{\begin{array}{ll}
\mathrm{div}(\tilde{\gamma}\nabla \tilde{u})=0&\text{in }\Omega\\
\tilde{u}=f &\text{on }\partial\Omega.
\end{array}
\right.
\end{equation}
By changing the Dirichlet datum $f$, one can perform two or more measurements. One often assumes that the
perturbation is well contained inside $\Omega$, that is, $\tilde{\gamma}$ coincides with the background conductivity $\gamma$ in a known neighbourhood of $\partial\Omega$, that is, outside $\tilde{\Omega}$, a known open set compactly contained in $\Omega$.

Since \cite{Mand}, it has been clear that one source of instability for elliptic inverse boundary value problems is due to the interior decay of solutions. We consider the solution $u$ to the Dirichlet problem in the unperturbed body, that is, with $\tilde{\gamma}$ replaced by the background conductivity $\gamma$, namely
\begin{equation}\label{dirpbmunpert}
\left\{\begin{array}{ll}
\mathrm{div}(\gamma\nabla u)=0&\text{in }\Omega\\
u=f &\text{on }\partial\Omega.
\end{array}
\right.
\end{equation}
The possibility to recover stably information on the unknown perturbation, using the additional measurement depending on the Dirichlet datum $f$, is directly related to the decay of $u$, or $\nabla u$, in the interior of the domain $\Omega$, in particular in $\tilde{\Omega}$, the region where the perturbation may be present. Therefore it is particularly important to establish decay properties of $u$ in $\tilde{\Omega}$, depending on the Dirichlet datum $f$ and the distance of $\tilde{\Omega}$ from $\partial\Omega$, which is exactly the issue we address in this paper.
 
 In \cite{Mand}, it was assumed that the domain $\Omega$ is $B_1$, a ball of radius $1$, and the background conductivity is homogeneous, $\gamma\equiv 1$. Then it was shown that the exponential instability of the inverse problem of Calder\'on is due to the fact that solutions to \eqref{dirpbmunpert} with boundary values $f$ given by spherical harmonics decay in the interior exponentially with respect to the degree of the spherical harmonic itself. We note that
a spherical harmonic is a Steklov eigenfunction for the Laplacian on $B_1$ with corresponding Steklov eigenvalue given by its degree. We also note that the degree is exactly equal to the frequency as defined in Definition~\ref{frequencydefin}.
The Weyl law on the asymptotic behaviour of Steklov eigenvalues is the other key ingredient. We recall that $\mu$ is a Steklov eigenvalue and $\phi|_{\partial\Omega}$ is its corresponding Steklov eigenfunction if
$\phi$ is a nontrivial solution to
\begin{equation}\label{Steklov}
\left\{\begin{array}{ll}
\mathrm{div}(\gamma\nabla \phi)=0&\text{in }\Omega\\
\gamma\nabla\phi\cdot\nu=\mu\phi &\text{on }\partial\Omega.
\end{array}
\right.
\end{equation}
We also recall that, for $f\in H^{1/2}(\partial\Omega)$, $f\neq 0$, we call its frequency the number
$$\mathrm{frequency}(f)=\frac{|f|^2_{H^{1/2}(\partial\Omega)}}{\|f\|^2_{L^2(\partial\Omega)}}.$$
We refer to Definition~\ref{frequencydefin} for a precise statement, here we just note that, if
$\phi|_{\partial\Omega}$ is the trace of a nontrivial solution to \eqref{Steklov}, then its frequency is essentially proportional to 
the Steklov eigenvalue $\mu$.

The ideas of \cite{Mand} have been generalised to other elliptic boundary value and scattering inverse problems in \cite{DC-R1,DC-R2} and to the parabolic case in \cite{DC-R-V}, showing that exponential instability unfortunately holds in all these cases.

Besides showing the instability nature of these problems, these results provide hints on the choice of optimal measurements, where optimality
may be in the sense of distinguishability as defined in \cite{Isaac}, see also \cite{G-I-N}. As suggested in these papers,
if one has at disposal a fixed and finite number $n$ of measurements, one should choose the first $n$ eigenfunctions of a suitable eigenvalue problem involving the perturbed conductivity $\tilde{\gamma}$ and the background one $\gamma$. Since the conductivity $\tilde{\gamma}$ is unknown, by using the arguments developed in \cite{Mand}
 the best choice should be to employ the first $n$ spherical harmonics, at least when the domain is a ball and the background conductivity is constant. In a general case, it seems reasonable to assume that the correct replacement for spherical harmonics is given by Steklov eigenfunctions. Indeed, the interior decay of the solution corresponding to a Steklov eigenfunction is very fast 
with respect to the Steklov eigenvalue, at least in a smooth case. For example, in \cite{His-Lut} it is shown that when $\partial\Omega$ is $C^{\infty}$ and $\gamma$ is $C^{\infty}$ the decay is faster than any power. Moreover, in the real-analytic case, the decay is still of exponential type, as shown first for surfaces in \cite{P-S-T} and then for higher dimensional manifolds in \cite{Gal-T}. 
Consequently, the information carried by the measurements corresponding to boundary data given by high order Steklov eigenfunctions rapidly degrades in the interior of the domain, thus it is of little help for the reconstruction of perturbations of the background conductivity far from the boundary.

However, from these examples it seems that the worst case scenario is when the domain and the background conductivity
are real-analyitc, because in this case the interior decay of the solution corresponding to a Steklov eigenfunction is indeed of exponential type with respect to the Steklov eigenvalue, or when the domain and the background conductivity are smooth, say $C^{\infty}$, since the interior decay is still very fast in this case.

Here, instead, we are interested in understanding the interior decay when the domain and the coefficients are not particularly smooth and also when $f$ is not a Steklov eigenfunction. In fact, in some occasions it might be very difficult to employ a Steklov eigenfunction and we wish to show that the decay might be actually due just to the frequency of the boundary datum $f$, without the much stronger assumption that $f$ is a Steklov eigenfunction. For example, this might be significant for the choice of optimal measurements in a partial data scenario, that is, when data are assigned and collected only on a given portion of the boundary.

We also wish to mention that, for the Calder\'on problem, the dependence of the distinguishability on the distance of $\tilde{\Omega}$ from the boundary of $\Omega$, rather than on the choice of boundary measurements, has been carefully analysed in \cite{Gar-Knu,A-S} in two dimensions and in \cite{Gar-Hyv} in higher dimensions. Our decay estimate, since the dependence on such a distance is explicitly given, may also be of interest in this kind of analysis of the instability.

We describe the main estimates we are able to prove, Theorems~\ref{mainthm} and \ref{mainthmbis}.
We assume that $\Omega$ is a $C^{1,1}$ domain and that $\gamma$ is Lipschitz continuous. We can also assume that $\gamma$ is a symmetric conductivity tensor, and not just a scalar conductivity, or that the underlying metric in $\Omega$ is not the Euclidean one but a Lipschitz Riemannian one. We call $\Phi$ the frequency of the Dirichlet boundary datum $f\in H^{1/2}(\partial\Omega)$, $f\neq 0$.
Whenever $f$ has zero mean on $\partial\Omega$, we may use another notion of frequency, which we call lower frequency and which is given by
$$\mathrm{lowefrequency}(f)=\frac{\|f\|^2_{L^2(\partial\Omega)}}{\|f\|^2_{H^{-1/2}(\partial\Omega)}}.$$
We refer to Definition~\ref{frequencydefinbis} for a precise statement. Here we point out that, if we call $\Phi_1$ the lower frequency of $f$, then $\Phi_1\leq \Phi$. On the other hand, if
$\phi|_{\partial\Omega}$ is the trace of a nontrivial solution to \eqref{Steklov} with $\mu>0$, then also its lower frequency is essentially proportional to 
the Steklov eigenvalue $\mu$.

For $d>0$ small enough, we call $\Omega^d$ the set
$$\Omega^d=\{x\in\Omega:\ \mathrm{dist}(x,\partial\Omega)>d\}.$$

The first result is the following.
We can find two positive constants $C_1$ and $C_2$, depending on $\Omega$, the Riemannian metric on it, and the coefficient $\gamma$, such that if $d \Phi\geq C_1$, then the function $u$ solving \eqref{dirpbmunpert} satisfies
\begin{equation}\label{mainestintro}
\int_{\Omega^d}\|\nabla u\|^2\leq C_2\frac{\int_{\Omega}\|\nabla u\|^2}{d\Phi}.
\end{equation}
We refer to Section~\ref{decaysec}, and in particular to Theorem~\ref{mainthm}, for the precise statement.

If we are interested in the decay of $u$ instead of its gradient, when $f$ has zero mean on $\partial\Omega$, we obtain an analogous result but we need to replace the frequency $\Phi$ with the lower frequency $\Phi_1$. Namely, we can find two positive constants $C_1$ and $C_2$, depending on $\Omega$, the Riemannian metric on it, and the coefficient $\gamma$, such that if $d \Phi_1\geq C_1$, then the function $u$ solving \eqref{dirpbmunpert} satisfies
\begin{equation}\label{mainestintrobis}
\int_{\partial\Omega^d} u^2\, d\sigma\leq C_2\frac{\int_{\partial\Omega}u^2\, d\sigma}{d\Phi_1}.
\end{equation}
We refer to Section~\ref{decaysec}, and in particular to Theorem~\ref{mainthmbis}, for the precise statement.

As an easy consequence of \eqref{mainestintrobis}, under the same assumptions we obtain that, for two positive constants $C_1$ and $C_2$, depending on $\Omega$, the Riemannian metric on it, and the coefficient $\gamma$, if $d \Phi_1\geq C_1$, then the function $u$ solving \eqref{dirpbmunpert} satisfies
\begin{equation}\label{mainestintroter}
\int_{\Omega^d}\|\nabla u\|^2\leq C_2\frac{\int_{\Omega}\|\nabla u\|^2}{d^2\Phi\Phi_1}.
\end{equation}
See Corollary~\ref{maincor} for the precise statement. We conclude that, if $f=\phi|_{\partial\Omega}$ is the trace of a nontrivial solution to \eqref{Steklov} with $\mu>0$, then, possibly with different constants $C_1$ and $C_2$, 
if $d \mu\geq C_1$, then the function $\phi$ solving \eqref{Steklov} satisfies
\begin{equation}\label{mainestintroquater}
\int_{\Omega^d}\|\nabla \phi\|^2\leq C_2\frac{\int_{\Omega}\|\nabla \phi\|^2}{d^2\mu^2},
\end{equation}
see Remark~\ref{Steklovremark} for a precise statement.

Let us briefly comment on the difference between these estimates.
Assuming that $f\in H^{1/2}(\partial\Omega)$, $f\neq 0$ with zero mean on $\partial\Omega$, since $\Phi_1\leq \Phi$, we have that $D$ in general decays faster than $H$. Actually, by Corollary~\ref{maincor}, we have that, as $\Phi_1$ grows, 
$D$ decays like $\Phi^{-1}\Phi_1^{-1}$, that is, at least like $\Phi_1^{-2}$, whereas $H$ decays like $\Phi_1^{-1}$. Moreover, if $f$ coincides with a Steklov eigenfunction with Steklov eigenvalue $\mu>0$, then, up to a constant, $\Phi$, $\Phi_1$ and $\mu$ are of the same order, therefore for Steklov eigenfunctions we obtain a decay of order $\mu^{-2}$, a result which is in accord with the estimate one can prove using the technique of \cite{His-Lut}.

In fact, an indication of the optimality of our decay estimates comes from the analysis developed in 
\cite{His-Lut} when $f=\phi|_{\partial\Omega}$ is a Steklov eigenfunction, with positive Steklov eigenvalue $\mu$. Following the idea of the proof of \cite[Theorem~1.1]{His-Lut}, it is evident that one can estimate $u(x)$, for any $x\in \Omega^d$, by a constant times $\mu^{-1}\|f\|_{H^{1/2}(\partial\Omega)}$ provided the Green's function $G_{\gamma}(x,\cdot)$ satisfies
\begin{equation}\label{Green}
\|\Lambda_{\gamma}(G_{\gamma}(x,\cdot))\|_{H^{1/2}(\partial\Omega)},\ \|\gamma\nabla G_{\gamma}(x,\cdot)\cdot\nu\|_{H^{1/2}(\partial\Omega)}\leq \tilde{C}
\end{equation}
where $\Lambda_{\gamma}$ is the so-called Dirichlet-to-Neumann map. 
Roughly speaking, \eqref{Green} corresponds to an $H^2$-bound of $G_{\gamma}(x,\cdot)$ away from $x$, which is what one obtains assuming the conductivity $\gamma$ is Lipschitz continuous.
Since $\|f\|_{H^{1/2}(\partial\Omega)}$ is of the order of $\sqrt{\mu}\|f\|_{L^2(\partial\Omega)}=\sqrt{\mu}\|u\|_{L^2(\partial\Omega)}$ and, as we already pointed out, the frequency and the lower frequency of $f$ are of the same order of $\mu$, one can obtain an estimate that is perfectly comparable with \eqref{mainestintrobis}. If one wishes to prove a decay of higher order, like $u(x)$ bounded by a constant times $\mu^{-2}\|f\|_{H^{1/2}(\partial\Omega)}$, by the same technique of \cite{His-Lut}, one should estimate the functions appearing in \eqref{Green} in terms of the $H^{3/2}(\partial\Omega)$ norm instead of the $H^{1/2}(\partial\Omega)$ norm, which corresponds to an $H^3$-bound of $G_{\gamma}(x,\cdot)$ away from $x$. Usually Lipschitz regularity of $\gamma$ is not enough to infer $H^3$-bounds, something like $C^{1,1}$ regularity would be required instead, 
therefore, under our weak regularity assumptions, our estimate \eqref{mainestintrobis} seems to be optimal even for Steklov eigenfunctions.

Another indication of the optimality of our decay estimates comes from the analysis developed in \cite{BLS}. In \cite{BLS} the authors introduce the so-called \emph{penetration function} and study its properties for two dimensional domains, in particular for the two dimensional unit ball. They are particularly interested in low regularity cases, thus they allow discontinuous conductivity tensors.
Their aim is to obtain estimates in homogenisation theory, but their results can be easily interpreted as distinguishability estimates with a finite number of boundary measurements for corresponding inverse boundary value problems.
In particular, using their notation, if $V_n$ is the space of trigonometric polynomials of degree $n$ on $\partial B_1(0)\subset \mathbb{R}^2$, $d\in (0,1)$ is a constant, and $A=\gamma$ is a symmetric conductivity tensor which is Lipschitz continuous, we can show that the penetration function $\Xi(V_n,d)$ satisfies, for a suitable constant $C$,
\begin{equation}\label{pfunction}
\Xi(V_n,d)\leq C(dn)^{-1}.
\end{equation}
In fact, for any $f$ which is orthogonal to $V_n$ in $L^2(\partial\Omega)$, we have that its frequency $\Phi$ is at least $n+1$ and also its lower frequency $\Phi_1$ is at least $n+1$. Therefore \eqref{pfunction} directly follows from \eqref{mainestintroter}.
Such a result considerably improves the estimate of \cite[Theorem~3.4]{BLS}, which is however valid for a wider class of conductivity tensors including discontinuous ones. Moreover, they give evidence by some explicit examples that,
when discontinuous conductivity tensors are allowed, a lower bound for the penetration function is of order $n^{-1/2}$.
It would be interesting to match such a lower bound by an estimate like \eqref{mainestintro} when $\gamma$ is discontinuous, but such an estimate would require a completely different method from the one used here.

About the technique we developed to obtain our estimates, let us begin by considering
\eqref{mainestintro}, where we use an ordinary differential equation argument that allows us to estimate the decay of
$$D(d)=\int_{\Omega^d}\gamma\|\nabla u\|^2$$
when $d$ is positive, and small enough. 
We closely follow the so-called frequency method introduced in \cite{Gar-Lin1}
to determine unique continuation properties of solutions to elliptic partial differential equations. In \cite{Gar-Lin1}, the local behaviour, near a point $x_0\in \Omega$, of a solution $u$ to $\mathrm{div}(\gamma\nabla u)=0$ in $\Omega$ was analysed, even in the case of a symmetric conductivity tensor $\gamma$.
A key point of the method was to reduce, locally near $x_0$, the elliptic equation with a symmetric conductivity tensor to an equation in a special Riemannian manifold with a scalar conductivity. By a special Riemannian manifold we mean one whose metric can be written in a special form in terms of polar coordinates centred at $x_0$. Such a reduction is made possible by the technique developed in \cite{AKS}.

Here we need to perform a similar construction, the only difference, and the main novelty, is that instead of considering a local modification near a point we consider a global one near the boundary of the domain. Indeed, in order to develop our analysis, we need that $\partial\Omega^d$ depends on $d$ smoothly enough or, equivalently, that the distance function from the boundary is smooth enough, say $C^{1,1}$, in a neighbourhood of the boundary. By \cite{Del-Zol}, see Theorem~\ref{Del-Zolthm}, this is true in the Euclidean setting provided $\partial\Omega$ is $C^{1,1}$ as well. In the Riemannian setting a similar result is much harder to prove. On the other hand, by exploiting the technique of \cite{AKS} and suitably changing the metric near the boundary, we can reduce to the case where the distance from the boundary, in the Riemannian metric, is smooth enough since it coincides with the distance from the boundary in the Euclidean metric in a neighbourhood of $\partial\Omega$.

We believe that
such a construction, besides being crucial for the proof of our decay estimates, is of independent interest and is one of the major achievement of the paper. The major part of the construction is contained in Proposition~\ref{viceversaprop}
and Theorem~\ref{AKSmethod}, with one interesting application developed in Proposition~\ref{normalderivativelemma}.

Our argument 
is based on the notion of frequency, which we essentially take from \cite{Gar-Lin1}, and which is given by
$$N(d)=\frac{D(d)}{H(d)}\qquad\text{where }H(d)=\int_{\partial\Omega^d}\gamma u^2\, d\sigma.$$
We note that $N(0)$ is of the same order of the frequency of the boundary datum $f$. We 
need to compute the derivative of $D$ and of $H$,
a task we perform following the analogous computations of \cite{Gar-Lin1}. In particular, for $D'(d)$ we use the coarea formula and a suitable version of the Rellich identity which is given in Lemma~\ref{Rellichid}. Instead, we compute $H'(d)$ by a straightforward application of Proposition~\ref{normalderivativelemma}.

The proof of \eqref{mainestintrobis} follows analogous lines of that of \eqref{mainestintro} by replacing $D$ with $H$ and $H$ with
$$E(d)=\int_{\Omega^d}\gamma u^2.$$
However there are some additional technical difficulties to be taken care of, see the proof of Theorem~\ref{mainthmbis} in Section~\ref{decaysec}. Moreover, the crucial link between the quotient $H(0)/E(0)$, which plays the role of $N(0)$, and the lower frequency $\Phi_1$ is provided by the estimate of Proposition~\ref{-1/2-2bound}.

The plan of the paper is as follows. In Section~\ref{prelsec} we present  the preliminary results that are needed for our analysis. In particular, we first discuss the regularity of domains and of the corresponding distance from the boundary, with the main result here being Theorem~\ref{Del-Zolthm} which is taken from \cite{Del-Zol}. We also give the precise definitions of frequencies we use. Then we review the Riemannian setting and the Dirichlet and Neumann problems for elliptic equations in the Euclidean and in the Riemannian setting, pointing out 
what happens if one suitably changes the underlying metric, see Remarks~\ref{conformalchange} and \ref{anisotropy}. For instance, Remark~\ref{anisotropy} allows us to pass from a symmetric conductivity tensor in the Euclidean setting to a scalar conductivity in the Riemannian one. We also briefly discuss Steklov eigenvalues and eigenfunctions.
In Section~\ref{distsec}, we investigate the distance function from the boundary in the Riemannian setting. Here the crucial result is Proposition~\ref{viceversaprop} which, together with Theorem~\ref{AKSmethod} and Remark~\ref{conformalchange}, allows us to assume, without loss of generality, that the distance function from the boundary in the Riemannian case has the same regularity as in the Euclidean case. Another important technical result in this section is Proposition~\ref{normalderivativelemma}.
Finally, in Section~\ref{decaysec}, we state and prove our main results, the decay estimates contained in
Theorems~\ref{mainthm} and \ref{mainthmbis} and Corollary~\ref{maincor}.

\subsubsection*{Acknowledgements} %
The authors are partly supported by GNAMPA, INdAM, through 2018 and 2019 projects. The authors wish to thank Eric Bonnetier for pointing them out reference \cite{BLS}.

\section{Preliminaries}\label{prelsec}

Throughout the paper the integer $N\geq 2$ will denote the space dimension. For any (column) vectors $v$, $w\in\mathbb{R}^N$, $\langle v,w\rangle = v^T w$ denotes the usual scalar product on $\mathbb{R}^N$. Here, and in the sequel,
for any matrix $A$, $A^T$ denotes its transpose.
For any $x=(x_1,\ldots,x_N)\in\mathbb{R}^N$,  we denote $x=(x',x_N)\in\mathbb{R}^{N-1}\times \mathbb{R}$.
We let $e_i$, $i=1,\ldots,N$, be the vectors of the canonical base and we call $\pi'$ the projection onto the first $(N-1)$ components and $\pi_N$ the projection onto the last one, namely, for any $x\in\mathbb{R}^N$,
$$\pi'(x)=x'=(x_1,\ldots,x_{N-1})\quad\text{and}\quad\pi_N(x)=x_N.$$

For any $s>0$ and any $x\in\mathbb{R}^N$, 
$B_s(x)$ denotes the open ball contained in $\mathbb{R}^N$ with radius $s$ and center $x$, whereas $B'_s(x')$ denotes the open ball contained in $\mathbb{R}^{N-1}$ with radius $s$ and center $x'$.
Finally, for any $E\subset \mathbb{R}^N$, we denote $B_s(E)=\bigcup_{x\in E}B_s(x)$.
For any Borel $E\subset\mathbb{R}^N$ we let $|E|=\mathcal{L}^N(E)$.
We call $\mathbb{M}^{N\times N}_{\mathrm{sym}}(\mathbb{R})$ the space of real-valued $N\times N$ symmetric matrices and by $I_N$ we denote the identity $N\times N$ matrix. We recall that we drop the dependence of any constant from the space dimension $N$.

\subsection{Regular domains and the distance from the boundary}

\begin{defin}
Let $\Omega\subset\mathbb{R}^N$ be a bounded open set. Let $k$ be a nonnegative integer and $0\leq\alpha\leq 1$.

We say that
$\Omega$ is of 
class $C^{k,\alpha}$
if for any $x\in\partial\Omega$ there exist a $C^{k,\alpha}$
function $\phi_x:\mathbb{R}^{N-1}\to\mathbb{R}$ and a neighbourhood $U_x$ of $x$
such that for any $y\in U_x$ we have, up to a rigid transformation depending on $x$,
$$y=(y',y_N)\in\Omega\quad \text{if and only if}\quad  y_N<\phi_x(y').$$

We also say that $\Omega$ is of 
class $C^{k,\alpha}$ with positive constants $r$ and $L$ if for any $x\in\partial\Omega$ we can choose $U_x=B_r(x)$ and $\phi_x$ such that $\|\phi_x\|_{C^{k,\alpha}(\mathbb{R}^{N-1})}\leq L$.
\end{defin}

\begin{oss}
If $\Omega\subset\mathbb{R}^N$, a bounded open set, is of 
class $C^{k,\alpha}$ then there exist positive constants $r$ and $L$ such that $\Omega$ is of 
class $C^{k,\alpha}$ with constants $r$ and $L$ with the further condition, when $k\geq 1$, that for any $x\in\partial\Omega$ we have $\nabla\phi_x (x')=0$.
\end{oss}

We note that a bounded open set of class $C^{0,1}$ is said to be of \emph{Lipschitz class} and that typically one assumes at least that $k+\alpha\geq 1$.

\begin{defin}\label{Eucldistfun}
Let $\Omega\subset\mathbb{R}^N$ be a bounded open set. 
For any $x\in \mathbb{R}^N$, its distance from the boundary of $\Omega$ is 
$$\mathrm{dist}(x,\partial\Omega)=\inf_{y\in\partial\Omega}\|x-y\|=\min_{y\in\partial\Omega}\|x-y\|.$$
We call $\varphi:\mathbb{R}^N\to\mathbb{R}$ the \emph{signed distance function} from the boundary of $\Omega$ as follows. For any $x\in \mathbb{R}^N$
$$\varphi(x)=\left\{\begin{array}{ll}\mathrm{dist}(x,\partial\Omega)
&\text{if }x\in\overline{\Omega},\\
-\mathrm{dist}(x,\partial\Omega)&\text{otherwise}.
\end{array}\right.$$
We call, for any $d\in\mathbb{R}$,
$$\Omega^d=\{x\in\mathbb{R}^N:\ \varphi(x)>d\}\quad\text{and}
\quad\partial\Omega^d=\{x\in\mathbb{R}^N:\ \varphi(x)=d\}.$$
Finally, for any $d>0$, we call
$$U^d=\{x\in\overline{\Omega}:\ \varphi(x)<d\}.$$
\end{defin}

The regularity of the signed distance function from the boundary has been thoroughly investigated in \cite{Del-Zol}.
Here we are interested in particular in the case of bounded open sets of class $C^{1,1}$ which is treated in \cite[Theorem~5.7]{Del-Zol}. Namely the following result holds true.

\begin{teo}\label{Del-Zolthm}
Let us fix positive constants $R$, $r$ and $L$.
Let $\Omega\subset B_R(0)\subset\mathbb{R}^N$ be a bounded open set of class $C^{1,1}$ with constants $r$ and $L$. Then there exists $\tilde{d}_0>0$, depending on $r$ and $L$ only, such that,
if we call $U=\{x\in \mathbb{R}^N:\ |\varphi(x)|<\tilde{d}_0\}$, for any $x\in U$ there exists a unique $y=P_{\partial\Omega}(x)\in \partial\Omega$ such that
$$\|x-P_{\partial\Omega}(x)\|=\mathrm{dist}(x,\partial\Omega).$$
Moreover, $\varphi$ is differentiable everywhere in $U$ and we have
\begin{equation}\label{varphiproper}
(\nabla\varphi(x))^T=-\nu(P_{\partial\Omega}(x))\text{ for any }x\in U,
\end{equation}
where $\nu$ denotes the exterior normal to $\Omega$, which we assume to be a column vector.
In particular,
$$\|\nabla\varphi\|=1\text{ in }U.$$

Finally, we have that $P_{\partial\Omega}\in C^{0,1}(U)$, with $C^{0,1}$ norm bounded by $r$, $L$ and $R$ only,
and, through \eqref{varphiproper},
we also have that  $\varphi\in C^{1,1}(U)$, with $C^{1,1}$ norm bounded by $r$, $L$ and $R$ only.
\end{teo}

\proof{.} It easily follows by using the arguments of the proof of \cite[Theorem~5.7]{Del-Zol}.\cvd

\smallskip

Let us note that, under the assumptions of Theorem~\ref{Del-Zolthm}, for any $0\leq |d|<\tilde{d}_0$,  we have that $\Omega^d$ is a bounded open set of class $C^{1,1}$ and
$\partial(\Omega^d)=\partial\Omega^d$. Moreover, for any $x\in \partial(\Omega^d)$,
if $\nu(x)$ denotes the exterior normal to $\Omega^d$, then
$$(\nabla\varphi(x))^T=-\nu(x)=-\nu(P_{\partial\Omega}(x)).$$

\begin{defin}\label{Omega_d}
Let $\Omega\subset\mathbb{R}^N$ be a bounded open set. We say that $A=A(x)\in \mathbb{M}^{N\times N}_{\mathrm{sym}}(\mathbb{R})$, $x\in\Omega$, is a \emph{symmetric tensor} in $\Omega$ if $A\in L^{\infty}(\Omega,\mathbb{M}^{N\times N}_{\mathrm{sym}}(\mathbb{R}))$.

We say that a symmetric tensor $A$ in $\Omega$ is \emph{Lipschitz} if $A\in C^{0,1}(\overline{\Omega},\mathbb{M}^{N\times N}_{\mathrm{sym}}(\mathbb{R}))$ and
that a symmetric tensor $A$ in $\Omega$ is \emph{uniformly elliptic} with constant 
$\lambda$, $0<\lambda<1$, if
$$\lambda\|\xi\|^2\leq \langle A(x)\xi,\xi \rangle\leq \lambda^{-1}\|\xi\|^2\quad\text{for almost any }x\in\Omega\text{ and any }\xi\in\mathbb{R}^N.$$
\end{defin}

If $\Omega$ is of class $C^{1,1}$ and $A$ is a Lipschitz conductivity tensor, we can extend $A$ outside $\Omega$ keeping it Lipschitz, and, in case, uniformly elliptic as well. Namely we have.

\begin{prop}\label{extensionprop}
Let us fix positive constants $R$, $r$ and $L$.
Let $\Omega\subset B_R(0)\subset\mathbb{R}^N$ be a bounded open set of class $C^{1,1}$ with constants $r$ and $L$. Let $A$ be a Lipschitz symmetric tensor in $\Omega$. Then there exists a Lipschitz symmetric tensor $\tilde{A}$ in $\mathbb{R}^N$ such that
$$\tilde{A}=A\text{ in }\Omega\quad\text{and}\quad \tilde{A}=I_N\text{ outside }B_{R+1}(0).$$
Moreover, the $C^{0,1}$ norm of $\tilde{A}$ on $\mathbb{R}^N$ depends on $r$, $L$, $R$ and the $C^{0,1}$ norm of $A$ on $\overline{\Omega}$. Finally, if $A$ is uniformly elliptic with constant $\lambda$, also $\tilde{A}$ is uniformly elliptic with the same constant $\lambda$.
\end{prop}

\proof{.} We sketch the idea of the construction. We pick $\tilde{d}_0$ and $U$ as in Theorem~\ref{Del-Zolthm} and we first extend $A$ in $\Omega\cup U$ as follows. We define, for any $x\in \Omega\cup U$,
$$\tilde{A}(x)=\left\{\begin{array}{ll}
A(x)
&\text{if }x\in\overline{\Omega}\\
A(P_{\partial\Omega}(x))&\text{if }x\in U\backslash\overline{\Omega}.
\end{array}
\right.
$$

Then we fix a cutoff function $\chi\in C^{\infty}(\mathbb{R})$ such that $\chi$ is increasing, $\chi(t)=0$ for any $t\leq -3\tilde{d}_0/4$ and $\chi(t)=1$ for any $t\geq 0$. We extend $\tilde{A}$ all over $\mathbb{R}^N$ as follows. We define, for any $x\in \mathbb{R}^N$,
$$\tilde{A}(x)=\chi(\varphi(x))\tilde{A}(x)+(1-\chi(\varphi(x)))I_N.$$
It is not difficult to check, with the help of Theorem~\ref{Del-Zolthm}, that such an extension satisfies the required properties.\cvd 

\smallskip

\subsection{Riemannian manifolds}

Let us consider the following definition of a Riemannian manifold $M$.

\begin{defin}\label{manifold}
 Let $\Omega\subset \mathbb{R}^N$ be a bounded open set of class $C^{1,1}$. Let  $G$ be a Lipschitz symmetric tensor in $\Omega$ which is uniformly elliptic with constant $\lambda$, $0<\lambda<1$.
For any $x\in\overline{\Omega}$, we denote as usual by $g_{i,j}(x)$ the elements of $G(x)$ and by $g^{i,j}(x)$ the elements of 
$G^{-1}(x)$, the inverse matrix of $G(x)$. Finally, we set $g(x)=|\det (G(x))|$.
We call $M$ the Riemannian manifold obtained by
endowing $\overline{\Omega}$ with the Lipschitz Riemannian metric whose tensor is given at any $x\in\overline{\Omega}$ by $g_{i,j}(x)dx_i\otimes dx_j$.

We finally say that $G$ is a scalar metric if $G=\theta I_N$ with $\theta \in C^{0,1}(\overline{\Omega})$, that is,
$g_{i,j}=\theta\delta_{i,j}$, where $\delta_{i,j}$ is the Kronecker delta.
\end{defin}

We recall the basic notation and properties of the Riemannian manifold $M$.
At any point $x\in\overline{\Omega}$, given any two (column) vectors $v$ and $w$, we denote
$$\langle v,w\rangle_M=\langle G(x)v,w \rangle$$
and, consequently,
$$\|v\|_M=\sqrt{\langle v,v\rangle_M}=\sqrt{\langle G(x)v,v\rangle}.$$
Clearly we have
$$
\sqrt{\lambda}\|v\|\leq\|v \|_M\leq \sqrt{\lambda^{-1}}\|v \|.
$$

For any $u\in L^1(\Omega)$, we have
$$\int_{\Omega}u(x)\, d_M(x)=\int_{\Omega}u(x)\sqrt{g(x)}\, dx.$$

If $h\in L^1(\partial\Omega)$, with respect to the surface measure $d\sigma$, that is, with respect to the $(N-1)$-dimensional Hausdorff measure, then
$$\int_{\partial \Omega}h(x)\, d\sigma_M(x)=\int_{\partial\Omega}h(x)\frac{\sqrt{g(x)}}{\alpha(x)}\, d\sigma(x),$$
where, for any $x\in\partial\Omega$,
$$\alpha(x)=\frac1{\sqrt{\langle G^{-1}(x)\nu(x),\nu(x)\rangle}},$$
$\nu(x)$ being the outer normal to the boundary. We call $\nu_M(x)=\alpha(x) G^{-1}(x)\nu(x)$, which is the outer normal to the boundary with respect to the Riemannian metric. In fact, $\|\nu_M(x)\|_M=1$ and
$\langle\tau,\nu_M(x)\rangle_M=0$ for any vector $\tau$ which is tangent to $\partial\Omega$ at the point $x$. 

At almost every $x\in\Omega$,
the intrinsic gradient of a function $u\in W^{1,1}(\Omega)$
is defined by 
$$\nabla_Mu(x)=\nabla u(x)G^{-1}(x)=g^{i,j}(x)\frac{\partial u}{\partial x_i}(x)e_j,$$
where we used the summation convention.
Let us note that, for any (column) vector $v$
$$\nabla u(x) v=\langle (\nabla u(x))^T,v\rangle=
\langle (\nabla_M u(x))^T,v\rangle_M.$$
Therefore,
\begin{multline}\label{gradcomp}
\|\nabla_M u(x)\|_M^2=\langle (\nabla_M u(x))^T,(\nabla_M u(x))^T\rangle_M\\=
\langle (\nabla u(x))^T,(\nabla_M u(x))^T\rangle=
\langle (\nabla u(x))^T,G^{-1}(x)(\nabla u(x))^T\rangle.
\end{multline}
Consequently,
\begin{equation}\label{normcomparison}
\sqrt{\lambda}\|\nabla u(x)\|\leq\|\nabla_M u(x)\|_M\leq \sqrt{\lambda^{-1}}\|\nabla u(x)\|.
\end{equation}

The intrinsic divergence of a vector field $X\in W^{1,1}(\Omega,\mathbb{R}^N)$ is defined, for almost every $x\in\Omega$, by
$$\mathrm{div}_M X(x)=\frac1{\sqrt{g(x)}}\mathrm{div}(\sqrt{g}X)(x).$$
For $X\in W^{1,1}(\Omega,\mathbb{R}^N)$, we have
$$\int_{\Omega}\mathrm{div}_M X(x)\, d_M(x)=\int_{\partial\Omega}\langle X(x),\nu_M(x)\rangle_M\, d\sigma_M(x).$$
Moreover, if $X\in W^{1,2}(\Omega,\mathbb{R}^N)$ and $\psi\in  W^{1,2}(\Omega)$, we have that
$$\mathrm{div}_M(X\psi)=\frac1{\sqrt{g}}\mathrm{div}(\sqrt{g}X)\psi+\nabla\psi X=
\mathrm{div}_M(X)\psi+\langle (\nabla_M \psi(x))^T,X\rangle_M.$$

Finally, the following version of the coarea formula holds true. Let $\varphi\in C^{1}(\overline{\Omega})$ be such that $\nabla \varphi\neq 0$ everywhere. Then for any $u\in L^1(\Omega)$, we have
$$\int_{\Omega}u(x)\, d_M(x)=\int_{\mathbb{R}}\left(\int_{\{x\in\Omega:\ \varphi(x)=t\}}\frac{u(x)}{\|\nabla_M\varphi(x)\|_M}\, d\sigma_M(x)\right)\, dt.$$

We call $\Gamma=\{\gamma:[0,1]\to \overline{\Omega}:\ \gamma\text{ is piecewise }C^1\}$.
For any curve $\gamma\in \Gamma$, we denote its Euclidean length as
$\mathrm{length}(\gamma)=\int_0^1\|\gamma'(t)\|\,dt$ and, analogously, its Riemannian length as
$$\mathrm{length}_M(\gamma)=\int_0^1\|\gamma'(t)\|_M\,dt.$$
We have that
$$\sqrt{\lambda}\,\mathrm{length}(\gamma)
\leq \mathrm{length}_M(\gamma)\leq \sqrt{\lambda^{-1}}\,\mathrm{length}(\gamma).$$
For any $x$ and $y\in\overline{\Omega}$, we call $\Gamma(x,y)=\{\gamma\in\Gamma:\ \gamma(0)=x\text{ and }\gamma(1)=y\}$ and define
$$d(x,y)=\inf_{\gamma\in \Gamma(x,y)}\, \mathrm{length}(\gamma)\quad\text{and}\quad d_M(x,y)=\inf_{\gamma\in \Gamma(x,y)}\, \mathrm{length}_M(\gamma).$$
Clearly
$$\sqrt{\lambda}\, d(x,y)
\leq d_M(x,y)\leq \sqrt{\lambda^{-1}}\, d(x,y),$$
whereas
\begin{equation}\label{lipd}
\|x-y\|\leq d(x,y)\leq C(\Omega)\|x-y\|,
\end{equation}
where $C(\Omega)$ is a constant depending on $\Omega$ only. If $\Omega$ satisfies the assumptions of Theorem~\ref{Del-Zolthm}, then $C(\Omega)$ depends on $r$, $L$ and $R$ only.

We finally define the distance from the boundary in the Riemannian case. Let 
$\varphi_M:\overline{\Omega}\to\mathbb{R}$ as follows. For any $x\in\overline{\Omega}$,
$$\varphi_M(x)=\mathrm{dist}_M(x,\partial\Omega)=\inf_{y\in\partial\Omega}d_M(x,y)=\min_{y\in\partial\Omega}d_M(x,y).$$
We observe that $\varphi$, the distance from the boundary in the Euclidean case that was defined in Definition~\ref{Eucldistfun}, satisfies
$$\varphi(x)=\mathrm{dist}(x,\partial\Omega)=\inf_{y\in\partial\Omega}d(x,y)=\min_{y\in\partial\Omega}d(x,y)\quad\text{for any }x\in\overline{\Omega}$$
and, consequently,
$$\sqrt{\lambda}\, \varphi(x)
\leq \varphi_M(x)\leq \sqrt{\lambda^{-1}}\, \varphi(x)\quad\text{for any }x\in\overline{\Omega}.$$

As in the Euclidean case, we adopt the following notation. For any $d\geq 0$, we define
$$\Omega_M^d=\{x\in\overline{\Omega}:\ \varphi_M(x)>d\}\quad\text{and}\quad
\partial\Omega_M^d=\{x\in\overline{\Omega}:\ \varphi_M(x)=d\}.$$
Moreover, when $d>0$, we call
$$U_M^d=\{x\in\overline{\Omega}:\ \varphi_M(x)<d\}.$$

We recall that Theorem~\ref{Del-Zolthm}, which easily follows from \cite[Theorem~5.7]{Del-Zol}, contains the regularity properties of $\varphi$, the (signed) distance function from the boundary in the Euclidean case. For the Riemannian metric, a corresponding regularity result for $\varphi_M$ is not easy to prove. We recall that fine regularity properties of the distance function from a general subset in a Riemannian manifold have been studied in \cite{Man-Men}.
In the next Section~\ref{distsec}, we study the properties of the distance function from the boundary in the Riemannian case.

\subsection{Definitions of frequencies of boundary data}

Let $\Omega\subset\mathbb{R}^N$ be a bounded Lipschitz domain. By domain we mean, as usual, an open and connected set.

We define the space of traces of $H^1(\Omega)$ functions on $\partial\Omega$ as 
$$H^{1/2}(\partial\Omega)=\{f=u|_{\partial\Omega}:\ u\in H^1(\Omega)\}.$$
We recall that $H^{1/2}(\partial\Omega)\subset L^2(\partial\Omega)$, with compact immersion.
By Poincar\'e inequality, an equivalent norm for $H^{1/2}(\partial\Omega)$, which we always adopt for simplicity, is given by the following
\begin{equation}\label{1/2norm}
\|f\|^2_{H^{1/2}(\partial\Omega)}=\|f\|^2_{L^2(\partial\Omega)}+|f|^2_{H^{1/2}(\partial\Omega)},
\end{equation}
where the seminorm is given by
\begin{equation}\label{1/2seminorm}
|f|^2_{H^{1/2}(\partial\Omega)}=\int_{\Omega}\|\nabla u_0(x)\|^2\, dx
\end{equation}
where $u_0\in H^1(\Omega)$ is the weak solution to
the following Dirichlet boundary value problem for the Laplace equation
\begin{equation}\label{Dirichlet0}
\left\{\begin{array}{ll}
\Delta u_0=0 & \text{in }\Omega\\
u_0=f & \text{on }\partial\Omega.
\end{array}\right.
\end{equation}

\begin{defin}\label{frequencydefin}
We call \emph{frequency} of a function $f\in H^{1/2}(\partial\Omega)$, with $f\neq 0$, the following quotient
\begin{equation}\label{frequency}
\mathrm{frequency}(f)=\frac{|f|^2_{H^{1/2}(\partial\Omega)}}{\|f\|^2_{L^2(\partial\Omega)}}\quad\text{for any }f\in H^{1/2}(\partial\Omega),\ f\neq 0.
\end{equation}
\end{defin}

\smallskip

We denote 
$L_{\ast}^2(\partial\Omega)=\{\psi\in L^2(\partial\Omega):\ \int_{\partial\Omega}\psi\, d\sigma=0\}$ and
$$H_{\ast}^{1/2}(\partial\Omega)=\left\{f\in H^{1/2}(\partial\Omega):\ \int_{\partial\Omega}f\, d\sigma=0\right\}.$$
We call $H^{-1/2}(\partial\Omega)$ the dual to $H^{1/2}(\partial\Omega)$ and
$$H_{\ast}^{-1/2}(\partial\Omega)=\{\eta\in H^{-1/2}(\partial\Omega):\ \langle \eta,1\rangle_{-1/2,1/2}=0\}.$$
By
$\langle\cdot,\cdot\rangle_{-1/2,1/2}$ we denote the duality between $H^{-1/2}(\partial\Omega)$ and $H^{1/2}(\partial\Omega)$.

By Poincar\'e inequality, we have that
$$\|f\|_{H_{\ast}^{1/2}(\partial\Omega)}=|f|_{H^{1/2}(\partial\Omega)}
\quad\text{for any }f\in H_{\ast}^{1/2}(\partial\Omega)$$
is an equivalent norm for $H_{\ast}^{1/2}(\partial\Omega)$ and, analogously,
$$
\|\eta\|_{H_{\ast}^{-1/2}(\partial\Omega)}=\sup_{\|\psi\|_{H^{1/2}_{\ast}(\partial\Omega)}=1}\langle\eta,\psi\rangle_{-1/2,1/2}
\quad
\text{for any }\eta\in H_{\ast}^{-1/2}(\partial\Omega)
$$
is an equivalent norm for $H_{\ast}^{-1/2}(\partial\Omega)$.

We observe that any $\eta\in L^2(\partial\Omega)$ is considered as an element of $H^{-1/2}(\partial\Omega)$ by setting
\begin{equation}\label{L2H-12}
\langle \eta,\psi\rangle_{-1/2,1/2}=\int_{\partial\Omega}\eta\psi\, d\sigma\quad\text{for any }\psi\in H^{1/2}(\partial\Omega).
\end{equation}
Moreover, if $\eta\in L_{\ast}^2(\partial\Omega)$ then $\eta\in H_{\ast}^{-1/2}(\partial\Omega)$. It is important to note that
here, and in the definitions of $L^2(\partial\Omega)$ and $L_{\ast}^2(\partial\Omega)$, we use the usual $(N-1)$-dimensional Hausdorff measure on $\partial\Omega$. In the sequel we adopt the same convention even if $\Omega$ is endowed with a Riemannian metric $G$ which is different from the Euclidean one. This simplifies the treatment of certain changes of variables for the Neumann problem or for the Steklov eigenvalue
problem, see Remark~\ref{Neumannpointwiseremark}.

\begin{defin}\label{frequencydefinbis}
We call \emph{lower frequency} of a function $f\in H^{1/2}_{\ast}(\partial\Omega)$, with $f\neq 0$, the following quotient
\begin{equation}\label{frequencybis}
\mathrm{lowfrequency}(f)=\frac{\|f\|^2_{L^2(\partial\Omega)}}{\|f\|^2_{H^{-1/2}_{\ast}(\partial\Omega)}}\quad\text{for any }f\in H^{1/2}_{\ast}(\partial\Omega),\ f\neq 0.
\end{equation}
Here
$$\|f\|_{H^{-1/2}_{\ast}(\partial\Omega)}=\sup_{\|\psi\|_{H^{1/2}_{\ast}(\partial\Omega)}=1}\int_{\partial\Omega}f\psi\, d\sigma
=\sup_{|\psi|_{H^{1/2}(\partial\Omega)}=1}\int_{\partial\Omega}f\psi\, d\sigma
.$$
\end{defin}

\smallskip

From this definition, we immediately infer that, for any $f\in H^{1/2}_{\ast}(\partial\Omega)$, with $f\neq 0$, we have
$$\|f\|^4_{L^2(\partial\Omega)}\leq \|f\|^2_{H^{-1/2}_{\ast}(\partial\Omega)}\|f\|^2_{H^{1/2}_{\ast}(\partial\Omega)}=
\|f\|^2_{H^{-1/2}_{\ast}(\partial\Omega)}|f|^2_{H^{1/2}(\partial\Omega)}
$$
hence
\begin{equation}\label{lowfrvsfr}
\mathrm{lowfrequency}(f)\leq \mathrm{frequency}(f)\quad\text{for any }f\in H^{1/2}_{\ast}(\partial\Omega),\ f\neq 0.
\end{equation}

\subsection{Boundary value problems for elliptic equations}

Let $\Omega\subset\mathbb{R}^N$ be a bounded Lipschitz domain.
We consider Dirichlet and Neumann problems in $\Omega$ for elliptic equations in divergence form, in the Euclidean and in the Riemannian setting.

Let $A=A(x)$ be a \emph{conductivity tensor} in $\Omega$, that is, $A$ is a symmetric tensor in $\Omega$ which is uniformly elliptic with some constant $\lambda_1$, $0<\lambda_1<1$.
If $A=\gamma I_N$, where $\gamma\in L^{\infty}(\Omega)$ satisfies
$$\lambda_1 \leq \gamma(x)\leq \lambda_1^{-1}\quad\text{for a.e. }x\in\Omega,$$
we say that $A$ (or $\gamma$) is a \emph{scalar conductivity}.

We say that a conductivity tensor $A$ is \emph{Lipschitz} if $A$ is a Lipschitz symmetric tensor. Analogously, 
$A$ (or $\gamma$) is a \emph{Lipschitz} scalar conductivity if $\gamma\in C^{0,1}(\overline{\Omega})$.

Let  $G$ be a Lipschitz symmetric tensor in $\Omega$ which is uniformly elliptic with constant $\lambda$, $0<\lambda<1$, and let $M$ be the corresponding Riemannian manifold on $\overline{\Omega}$ as in Definition~\ref{manifold}.

In this subsection we adopt the following assumption.

\begin{assum}\label{scalarassum}
We assume that either $A$ is a scalar conductivity tensor, that is, $A=\gamma I_N$ with $\gamma\in L^{\infty}(\Omega)$, or
$G$ is a scalar metric, that is, $G=\theta I_N$ with $\theta\in C^{0,1}(\overline{\Omega})$.
\end{assum}

For any $f\in H^{1/2}(\partial\Omega)$, let $u\in H^1(\Omega)$ be the weak solution to the Dirichlet boundary value problem
\begin{equation}\label{Dirichlet}
\left\{\begin{array}{ll}
\mathrm{div}_M(A\nabla_M u)=0 & \text{in }\Omega\\
u=f & \text{on }\partial\Omega.
\end{array}\right.
\end{equation}
We recall that $u\in H^1(\Omega)$ solves \eqref{Dirichlet} if $u=f$ on $\partial\Omega$ in the trace sense and
$$\int_{\Omega}\langle A(x)(\nabla_M u(x))^T,(\nabla_M \psi(x))^T\rangle_M\, d_M(x)=0\quad\text{for any }\psi\in H^1_0(\Omega).$$
For the sake of simplicity, we sometimes drop the transpose in the sequel, considering, with a small abuse of notation, the gradient as a column vector.

The following remark holds true.

\begin{oss}\label{freqrem}
Let $u$ and $u_0$ be the solution to \eqref{Dirichlet} and \eqref{Dirichlet0}, respectively.
Then there exists a constant $c_1$, $0<c_1<1$ depending on $\lambda$ and $\lambda_1$ only, such that
\begin{equation}\label{frequencyvsintegral}
c_1
\int_{\Omega}\|\nabla u_0(x)\|^2\, dx\leq \int_{\Omega}\langle A(x)\nabla_M u(x),\nabla_M u(x)\rangle_M\, d_M(x)\leq 
c_1^{-1}
\int_{\Omega}\|\nabla u_0(x)\|^2\, dx.
\end{equation}
In fact, on the one hand, by the Dirichlet principle,
$$\int_{\Omega}\langle A(x)\nabla_M u(x),\nabla_M u(x)\rangle_M\, d_M(x)\geq c_1
\int_{\Omega}\|\nabla u(x)\|^2\, dx\geq c_1 \int_{\Omega}\|\nabla u_0(x)\|^2\, dx.$$
On the other hand, correspondingly we have
\begin{multline*}
\int_{\Omega}\langle A(x)\nabla_M u(x),\nabla_M u(x)\rangle_M\, d_M(x)\leq
\int_{\Omega}\langle A(x)\nabla_M u_0(x),\nabla_M u_0(x)\rangle_M\, d_M(x)\\\leq
c_1^{-1}\int_{\Omega}\|\nabla u_0(x)\|^2\, dx.
\end{multline*}
\end{oss}

As a consequence of Remark~\ref{freqrem}, we can define equivalent $H^{1/2}(\partial\Omega)$ norm and seminorm which are given by, for any $f\in H^{1/2}(\partial\Omega)$,
\begin{equation}
|f|^2_{H^{1/2}_A(\partial\Omega)}=\int_{\Omega}\langle A(x)\nabla_M u(x),\nabla_M u(x)\rangle_M\, d_M(x),
\end{equation}
where $u$ solves \eqref{Dirichlet}, and
\begin{equation}
\|f\|^2_{H^{1/2}_A(\partial\Omega)}=\|f\|^2_{L^2(\partial\Omega)}+|f|^2_{H^{1/2}_A(\partial\Omega)}.
\end{equation}
We can also define an equivalent $H^{-1/2}(\partial\Omega)$ norm given by, for any $\eta\in
H^{-1/2}(\partial\Omega)$,
$$\|\eta\|_{H^{-1/2}_A(\partial\Omega)}=\sup_{\|\psi\|_{H^{1/2}_A(\partial\Omega)}=1}\langle \eta,\psi\rangle_{-1/2,1/2}.$$
We note that here we drop the dependence on the metric $M$, although the seminorm, and thus the norms as well, clearly also depends on it.

Analogously,
$$\|f\|_{H_{\ast,A}^{1/2}(\partial\Omega)}=|f|_{H^{1/2}_A(\partial\Omega)}
\quad\text{for any }f\in H_{\ast}^{1/2}(\partial\Omega)$$
is an equivalent norm for $H_{\ast}^{1/2}(\partial\Omega)$ and
$$
\|\eta\|_{H_{\ast,A}^{-1/2}(\partial\Omega)}=\sup_{\|\psi\|_{H^{1/2}_{\ast,A}(\partial\Omega)}=1}\langle\eta,\psi\rangle_{-1/2,1/2}
\quad\text{for any }\eta\in H_{\ast}^{-1/2}(\partial\Omega)
$$
is an equivalent norm for $H_{\ast}^{-1/2}(\partial\Omega)$.

For any $\eta\in H_{\ast}^{-1/2}(\partial\Omega)$, let $v\in H^1(\Omega)$ be the solution to the Neumann boundary value problem 
\begin{equation}\label{Neumann}
\left\{\begin{array}{ll}
\mathrm{div}_M(A\nabla_M v)=0 & \text{in }\Omega\\
\langle A\nabla_M v,\nu_M\rangle_M=\eta & \text{on }\partial\Omega\\
\int_{\partial\Omega}v\, d\sigma=0. &
\end{array}\right.
\end{equation}
By a solution we mean $v\in H^1(\Omega)$ such that
$v|_{\partial\Omega}\in H_{\ast}^{1/2}(\partial\Omega)$ and that
$$\int_{\Omega} \langle A\nabla_M v,\nabla_M \psi\rangle_M=\langle \eta,\psi|_{\partial\Omega}\rangle_{-1/2,1/2}\quad\text{for any }\psi\in H^1(\Omega).$$
We also note that, for simplicity and by a slight abuse of notation, we denote
$A u_{\nu_M}=\langle A\nabla_M v,\nu_M\rangle_M$. Such a notation is actually correct when $A=\gamma I_N$ is a scalar conductivity. In fact, in this case,
$$A u_{\nu_M}=\gamma u_{\nu_M}$$
where $u_{\nu_M}$ is the (exterior) normal derivative of $u$ with respect to $\Omega$ which, in the Riemannian setting, is given by
$$u_{\nu_M}=\langle(\nabla_Mu)^T,\nu_M\rangle_M=\nabla u \nu_M.$$

By Poincar\'e inequality and Lax-Milgram lemma, we have that there exists a unique solution both to \eqref{Dirichlet} and to \eqref{Neumann}. Moreover,
there exists a constant $c_2$, $0<c_2<1$ depending on $\Omega$, $\lambda$ and $\lambda_1$ only, such that
for any $f\in H^{1/2}(\partial\Omega)$
$$c_2\|f\|_{H^{1/2}(\partial\Omega)}\leq \|u\|_{H^1(\Omega)}\leq c_2^{-1}\|f\|_{H^{1/2}(\partial\Omega)}$$
and for any $\eta\in H_{\ast}^{-1/2}(\partial\Omega)$
$$c_2\|\eta\|_{H^{-1/2}(\partial\Omega)}\leq \|v\|_{H^1(\Omega)}\leq c_2^{-1}\|\eta\|_{H^{-1/2}(\partial\Omega)}.$$
If $\Omega\subset B_R(0)$ is Lipschitz with positive constants $r$ and $L$, the dependence of $c_2$ on $\Omega$ is just through the constants $r$, $L$ and $R$.

Let $\Lambda:H^{1/2}(\partial\Omega)\to H_{\ast}^{-1/2}(\partial\Omega)$ be the linear operator such that
$$\Lambda(f)=\langle A\nabla_M v,\nu_M\rangle_M\quad \text{for any }f\in H^{1/2}(\partial\Omega)$$
where $u$ solves \eqref{Dirichlet}. Here, we mean
$$\langle \langle A\nabla_M v,\nu_M\rangle_M,\tilde{\psi}\rangle_{-1/2,1/2}=\int_{\Omega}\langle A\nabla_Mu,\nabla_M\psi\rangle_M\, d_M\quad\text{for any }\tilde{\psi}\in H^{1/2}(\partial\Omega),$$
where $\psi$ is any $H^1(\Omega)$ function such that $\psi|_{\partial\Omega}=\tilde{\psi}$.
 We infer that $\Lambda$ restricted to $H_{\ast}^{1/2}(\partial\Omega)$ is invertible and both $\Lambda$ and $\Lambda^{-1}:H_{\ast}^{-1/2}(\partial\Omega)\to H_{\ast}^{1/2}(\partial\Omega)$ are bounded operators with norms bounded by constants depending on $\Omega$, $\lambda$ and $\lambda_1$ only. As usual we refer to $\Lambda$ as the Dirichlet-to-Neumann map and to $\Lambda^{-1}$ as the Neumann-to-Dirichlet map.

We are interested in eigenvalues and eigenfunctions of the Dirichlet-to-Neumann map $\Lambda$, which coincides with the so-called Steklov eigenvalues and eigenfunctions. Namely, we say that $\mu\in \mathbb{C}$ and $\phi\in L^2(\partial\Omega)$, with $\phi\neq 0$ are, respectively, a Steklov eigenvalue and its corresponding eigenfunction if
there exists $w\in H^1(\Omega)$ such that $w=\phi$ on $\partial\Omega$ and $w$ satisfies
\begin{equation}\label{Steklovequation}
\left\{
\begin{array}{ll}
\mathrm{div}_M(A\nabla_M w)= 0 & \text{in }\Omega\\
\langle A\nabla_M v,\nu_M\rangle_M
=\mu w &\text{on }\partial\Omega,
\end{array}\right.
\end{equation}
that is,
$$\int_{\Omega} \langle A\nabla_M w,\nabla_M \psi\rangle_M\, d_M=\langle \mu w|_{\partial\Omega},\psi|_{\partial\Omega}\rangle_{-1/2,1/2}=\int_{\partial\Omega}\mu w\psi\, d\sigma
\quad\text{for any }\psi\in H^1(\Omega).$$
In other words, $\phi$ satisfies $\Lambda(\phi)=\mu \phi$.
Clearly \eqref{Steklovequation} is satisfied by $\mu=0$ and $w$ a constant function. It is well-known that the Steklov eigenvalues form
an increasing sequence of real numbers
$$0=\mu_0<\mu_1\leq \mu_2\leq \ldots \leq \mu_n\leq \ldots$$
such that $\lim_n\mu_n=+\infty$. For any $n\geq 0$, we can find a corresponding eigenfunction $\phi_n$, normalised in such a way that $\|\phi_n\|_{L^2(\partial\Omega)}=1$, such that
$\{\phi_n\}_{n\geq 0}$ is an orthonormal basis of $L^2(\partial\Omega)$ and
$\{\phi_n\}_{n\in\mathbb{N}}$ is an orthonormal basis of $L^2_{\ast}(\partial\Omega)$. Moreover,
$\{\phi_n/\sqrt{1+\mu_n}\}_{n\geq 0}$ and $\{\phi_n/\sqrt{1+\mu_n}\}_{n\in\mathbb{N}}$ are an orthonormal basis of $H^{1/2}(\partial\Omega)$ and $H^{1/2}_{\ast}(\partial\Omega)$, respectively, 
with respect to the $H^{1/2}_A(\partial\Omega)$ norm.
Finally, we call $\{\psi_n=\phi_n/\sqrt{\mu_n}\}_{n\in\mathbb{N}}$ and we note that it is an 
orthonormal basis of $H^{1/2}_{\ast}(\partial\Omega)$
with respect to the $H^{1/2}_{\ast,A}(\partial\Omega)$ norm.

If $\phi\in H_{\ast}^{1/2}(\partial\Omega)$ is a Steklov eigenfunction with eigenvalue $\mu$, and $w$ is the corresponding solution to \eqref{Steklovequation}, then
$$\mu=
\frac{\int_{\partial\Omega}\mu\phi^2\, d\sigma}{\int_{\partial\Omega}\phi^2\, d\sigma}=
\frac{\int_{\Omega}\langle A\nabla_Mw,\nabla_Mw\rangle_M\, d_M}{\int_{\partial\Omega}\phi^2\, d\sigma},$$
hence by Remark~\ref{freqrem} we have, with the same constant $c_1$,
\begin{equation}\label{Stekvsfreq}
c_1\,\mathrm{frequency}(\phi)\leq \mu \leq c_1^{-1}\,\mathrm{frequency}(\phi).
\end{equation}
An important property of Steklov eigenfunctions is that their frequency and lower frequency are of the same order. In fact, for $\mu>0$ we have $\phi\in H^{1/2}_{\ast}(\partial\Omega)$ and, setting $\int_{\partial\Omega}\phi^2=1$,
$$c_1|\phi|^2_{H^{1/2}(\partial\Omega)}\leq\|\phi\|^2_{H^{1/2}_{\ast,A}(\partial\Omega)}=\mu\leq c_1^{-1}|\phi|^2_{H^{1/2}(\partial\Omega)},$$
therefore
$$c_1\|\phi\|^2_{H^{-1/2}_{\ast}(\partial\Omega)}\leq\|\phi\|^2_{H^{-1/2}_{\ast,A}(\partial\Omega)}=\mu^{-1}\leq c_1^{-1}\|\phi\|^2_{H^{-1/2}_{\ast}(\partial\Omega)},$$
and, finally,
\begin{equation}\label{Stekvslowfreq}
c_1\,\mathrm{lowfrequency}(\phi)\leq \mu \leq c_1^{-1}\,\mathrm{lowfrequency}(\phi).
\end{equation}

Although their proofs are elementary, and actually quite similar, the next two remarks are crucial.

\begin{oss}\label{conformalchange}
Let $A$ be a conductivity tensor in $\Omega$ which is uniformly elliptic with some constant $\lambda_1$, $0<\lambda_1<1$.
Let  $G$ be a Lipschitz symmetric tensor in $\Omega$ which is uniformly elliptic with constant $\lambda$, $0<\lambda<1$, and let $M$ be the corresponding Riemannian manifold on $\overline{\Omega}$.
Let Assumption~\ref{scalarassum} be satisfied.

Let us take $\eta_1\in C^{0,1}(\overline{\Omega})$ such that $\lambda_1\leq \eta_1\leq\lambda^{-1}_1$ in $\overline{\Omega}$,
for some constant $\lambda_1$, $0<\lambda_1<1$.
Let us define $\tilde{G}=\eta_1 G$ and let us consider the Riemannian manifold $\tilde{M}$ obtained by endowing $\overline{\Omega}$ with the Lipschitz Riemannian metric given by $\tilde{G}$.

We define
$$\tilde{A}=\eta_1^{(2-N)/2}A,$$
and we note that $\tilde{A}=A$ if $N=2$.

Then, for any $\psi_1$, $\psi_2\in H^1(\Omega)$ we have
$$\int_{\Omega}\langle A \nabla_M\psi_1,\nabla_M\psi_2\rangle_M\, d_M=
\int_{\Omega}\langle \tilde{A} \nabla_{\tilde{M}}\psi_1,\nabla_{\tilde{M}}\psi_2\rangle_{\tilde{M}}\, d_{\tilde{M}}.$$
\end{oss}

The next remark shows that, under Assumption~\ref{scalarassum} and if $A$ is Lipschitz, we can always assume that
the conductivity tensor is a scalar conductivity, up to changing the Riemannian metric. For example, this applies when $A$ is a Lipschitz conductivity tensor and the metric is the Euclidean one. Namely we have the following.

\begin{oss}\label{anisotropy}
Let $A$ be a Lipschitz conductivity tensor in $\Omega$ which is uniformly elliptic with some constant $\lambda_1$, $0<\lambda_1<1$.
Let  $G$ be a Lipschitz symmetric tensor in $\Omega$ which is uniformly elliptic with constant $\lambda$, $0<\lambda<1$, and let $M$ be the corresponding Riemannian manifold on $\overline{\Omega}$.
Let Assumption~\ref{scalarassum} be satisfied.

We call $A_1=\sqrt{g}A G^{-1}$ and $\gamma_1=(\det{A_1})^{1/N}$
so that $A_1=\gamma_1\hat{A}_1$ with $\det{\hat{A}_1}\equiv 1$.

If $N>2$, we define $\tilde{A}\equiv I_N$  and
$$\tilde{G}=(\det(A_1))^{1/(N-2)}A_1^{-1}.$$

If $N=2$, we define $\tilde{A}\equiv \gamma_1 I_N$ and
$$\tilde{G}=\hat{A}_1^{-1}.$$

Let us consider the Riemannian manifold $\tilde{M}$ obtained by endowing $\overline{\Omega}$ with the Lipschitz Riemannian metric given by $\tilde{G}$.

Then, for any $\psi_1$, $\psi_2\in H^1(\Omega)$ we have
$$\int_{\Omega}\langle A \nabla_M\psi_1,\nabla_M\psi_2\rangle_M\, d_M=
\int_{\Omega}\langle \tilde{A} \nabla_{\tilde{M}}\psi_1,\nabla_{\tilde{M}}\psi_2\rangle_{\tilde{M}}\, d_{\tilde{M}}.$$
\end{oss}

Both for the case of Remark~\ref{conformalchange} and the one of Remark~\ref{anisotropy}, we infer the following consequences.

Fixed $f\in H^{1/2}(\partial\Omega)$, let $u$ be the solution to \eqref{Dirichlet}. Then $u$ solves
\begin{equation}\label{Dirichletmod}
\left\{\begin{array}{ll}
\mathrm{div}_{\tilde{M}}(\tilde{A}\nabla_{\tilde{M}} u)=0 & \text{in }\Omega\\
u=f & \text{on }\partial\Omega.
\end{array}\right.
\end{equation}

Analogously, fixed $\eta\in H_{\ast}^{-1/2}(\partial\Omega)$,
let $v$ be the solution to \eqref{Neumann}. Then $v$ solves
\begin{equation}\label{Neumannmod}
\left\{\begin{array}{ll}
\mathrm{div}_{\tilde{M}}(\tilde{A}\nabla_{\tilde{M}} v)=0 & \text{in }\Omega\\
\langle \tilde{A}\nabla_{\tilde{M}} v,\nu_{\tilde{M}}\rangle_{\tilde{M}}=\eta & \text{on }\partial\Omega\\
\int_{\partial\Omega}v\, d\sigma=0. &
\end{array}\right.
\end{equation}

Finally, if $w$ solves \eqref{Steklovequation} for a constant $\mu$, then
 $w$ solves
\begin{equation}\label{Steklovmod}
\left\{\begin{array}{ll}
\mathrm{div}_{\tilde{M}}(\tilde{A}\nabla_{\tilde{M}} w)=0 & \text{in }\Omega\\
\langle \tilde{A}\nabla_{\tilde{M}} w,\nu_{\tilde{M}}\rangle_{\tilde{M}}=\mu w & \text{on }\partial\Omega.
\end{array}\right.
\end{equation}

We conclude this section by investigating the regularity of the solutions to \eqref{Dirichlet}, \eqref{Neumann} and \eqref{Steklovequation}. We need stronger assumptions on the domain $\Omega$ and the conductivity tensor $A$. Namely we assume the following till the end of the section.

Let us fix positive constants $R$, $r$, $L$, $C_0$, $C_1$, $\lambda$ and $\lambda_1$, 
with $0<\lambda<1$ and $0<\lambda_1<1$. We refer to these constants as the
\emph{a priori data}.

Let $\Omega\subset B_R(0)\subset\mathbb{R}^N$ be a bounded domain of class $C^{1,1}$ with constants $r$ and $L$.

Let  $G$ be a Lipschitz symmetric tensor in $\Omega$ which is uniformly elliptic with constant $\lambda$
and such that $\|G\|_{C^{0,1}(\overline{\Omega})}\leq C_0$.

Let $A$ be a Lipschtitz conductivity tensor in $\Omega$ 
which is uniformly elliptic with constant $\lambda_1$
and such that
$\|A\|_{C^{0,1}(\overline{\Omega})}\leq C_1$.

We suppose that Assumption~\ref{scalarassum} holds.
We note that, without loss of generality, through Remark~\ref{anisotropy}, we could just assume that $A$ is a scalar conductivity.

The first remark is that, by standard regularity estimates for elliptic equations,
if $u$ is any weak solution to $\mathrm{div}_M(A\nabla_M u)=0$ in $\Omega$, then
 $u\in H^2_{\mathrm{loc}}(\Omega)$ and the equation is satisfied pointwise almost everywhere in $\Omega$.

Here we are interested on the conditions that guarantee that our solutions are actually belonging to $H^2(\Omega)$.

We adopt the standard definition of $H^{3/2}(\partial\Omega)$, see for example \cite{G}, and by $H_{\ast}^{3/2}(\partial\Omega)$ we denote the elements of $H^{3/2}(\partial\Omega)$ with zero mean on $\partial\Omega$.
Let $u$ be the solution to \eqref{Dirichlet} with boundary datum $f\in H^{1/2}(\partial\Omega)$ and
$v$ the solution to \eqref{Neumann} with boundary datum $\eta\in H^{-1/2}_{\ast}(\partial\Omega)$. The following regularity properties hold true. 

\begin{prop}\label{regprop}
There exist a positive constants $c_3$, $0<c_3<1$ depending on the a priori data only,
such that
for any $f\in H^{3/2}(\partial\Omega)$
\begin{equation}\label{aineq}
c_3\|f\|_{H^{3/2}(\partial\Omega)}\leq \|u\|_{H^2(\Omega)}\leq c_3^{-1}\|f\|_{H^{3/2}(\partial\Omega)}
\end{equation}
and for any $\eta\in H_{\ast}^{1/2}(\partial\Omega)$
\begin{equation}\label{bineq}
c_3\|\eta\|_{H^{1/2}(\partial\Omega)}\leq \|v\|_{H^2(\Omega)}\leq c_3^{-1}\|\eta\|_{H^{1/2}(\partial\Omega)}.
\end{equation}
In \eqref{bineq}, we can replace $\|\eta\|_{H^{1/2}(\partial\Omega)}$ with $\|\eta\|_{H^{1/2}_{\ast}(\partial\Omega)}$,
 $\|\eta\|_{H^{1/2}_A(\partial\Omega)}$ or $\|\eta\|_{H^{1/2}_{\ast,A}(\partial\Omega)}$.

As a consequence, $\Lambda$ is bounded between $H_{\ast}^{3/2}(\partial\Omega)$ and $H_{\ast}^{1/2}(\partial\Omega)$, with a bounded inverse, and their norms are bounded by constants depending on the a priori data only.
\end{prop}

Before sketching the proof of this standard regularity result, we state the following important remark.

\begin{oss}\label{Neumannpointwiseremark}
Let $v\in H^2(\Omega)$ be a solution to $\mathrm{div}_M(A\nabla_M v)=0$ in $\Omega$. Then $\nabla v\in H^1(\Omega)$, therefore $\nabla v$ is well-defined, in the trace sense, on $\partial\Omega$.
It follows that $A v_{\nu_M}$ is well-defined for instance in $L^2(\partial\Omega)$.
 Moreover, using integration by parts, we conclude that for any $\psi\in H^1(\Omega)$ we have
$$\int_{\Omega} \langle A\nabla_M v,\nabla_M \psi\rangle_M\, d_M=
\int_{\partial\Omega}A v_{\nu_M}\psi\, d\sigma_M=\int_{\partial\Omega}\eta\psi\, d\sigma$$
where
\begin{equation}\label{Neumannpointwise}
\eta= \frac{\sqrt{g}}{\alpha}A v_{\nu_M}=
\sqrt{g}\langle A \nabla_M v,\nu\rangle.
\end{equation}
Therefore, in the Riemannian setting, the Neumann condition
$$A v_{\nu_M}=\eta\quad\text{on }\partial\Omega$$
is in general not valid in a pointwise or $L^2$ sense, even when both $Av_{\nu_M}$ and $\eta$ are well-defined as $L^2(\partial\Omega)$ functions. The correct pointwise or $L^2$ boundary condition is given in \eqref{Neumannpointwise}.
\end{oss}

\proof{ of Proposition~\textnormal{\ref{regprop}}.} This result is essentially proved in \cite{G}.

Using for instance \cite[Theorem~1.5.1.2]{G} and \cite[Theorem~1.5.1.3]{G}, with the help of Remark~\ref{Neumannpointwiseremark}, we immediately infer that the
left inequalities of \eqref{aineq} and \eqref{bineq} hold true.

The right inequalities of \eqref{aineq} and \eqref{bineq} easily follow by \cite[Corollary~2.2.2.4]{G} and \cite[Corollary~2.2.2.6]{G}.\cvd

\smallskip

An important consequence of Proposition~\ref{regprop} for Steklov eigenfunctions is the following.

\begin{cor}\label{Stekest}
Let $\phi\in H_{\ast}^{1/2}(\partial\Omega)$ be a Steklov eigenfunction with eigenvalue $\mu>0$ and let
$w$ be the corresponding solution to \eqref{Steklovequation}.

Then
\begin{equation}\label{H2steklov}
c_3^2\mu^2(1+c_1\mu)\|\phi\|^2_{L^2(\partial\Omega)}\leq
\|w\|^2_{H^2(\Omega)}\leq c_3^{-2}\mu^2(1+c_1^{-1}\mu)\|\phi\|^2_{L^2(\partial\Omega)},
\end{equation}
where $c_1$ is as in \eqref{Stekvsfreq} and $c_3$ is as in Proposition~\textnormal{\ref{regprop}}, thus they depend on the a priori data only.
\end{cor}

\proof{.} By \eqref{Stekvsfreq}, we have that
$$c_1\mu\|\phi\|^2_{L^2(\partial\Omega)} \leq|\phi|^2_{H^{1/2}(\partial\Omega)}\leq c_1^{-1}\mu\|\phi\|^2_{L^2(\partial\Omega)}.$$
Therefore
$$(1+c_1\mu)\|\phi\|^2_{L^2(\partial\Omega)}\leq
\|\phi\|^2_{H^{1/2}(\partial\Omega)}\leq (1+c_1^{-1}\mu)\|\phi\|^2_{L^2(\partial\Omega)}.$$
Then the result follows by 
Proposition~\ref{regprop}, in particular by \eqref{bineq} with $\eta=\mu \phi$.\cvd

\smallskip

Finally, we state and prove the following result.

\begin{prop}\label{-1/2-2bound}
There exists a constant $C_2$, depending on the a priori data only, such that for any $f\in H^{1/2}(\partial\Omega)$
we have
\begin{equation}\label{-1/2-2boundest}
\|u\|_{L^2(\Omega)}\leq C_2\|f\|_{H^{-1/2}(\partial\Omega)},
\end{equation}
where $u$ is the solution to \eqref{Dirichlet}.
\end{prop}

\proof{.} Without loss of generality, we can restrict our attention to $f\in H_{\ast}^{1/2}(\partial\Omega)$ and
we can replace the $H^{-1/2}(\partial\Omega)$ norm with the
$H^{-1/2}_{\ast,A}(\partial\Omega)$ norm.
Given $f\in H^{1/2}_{\ast}(\partial\Omega)$,
we can find a sequence $\{\alpha_n\}_{n\in\mathbb{N}}$ of real numbers such that
$$f=\sum_{n\in\mathbb{N}}\alpha_n\psi_n\quad\text{and}\quad\|f\|^2_{H^{1/2}_{\ast,A}(\partial\Omega)}=\sum_{n\in\mathbb{N}}\alpha_n^2.$$
Furthermore, it is easy to infer that
$$\|f\|^2_{H^{-1/2}_{\ast,A}(\partial\Omega)}=\sum_{n\in\mathbb{N}}\frac{\alpha_n^2}{\mu_n^2}.$$

We have that
$$\Lambda(f)=\sum_{n\in\mathbb{N}}\alpha_n\mu_n\psi_n,$$
therefore
$$\|\Lambda(f)\|^2_{H^{1/2}_{\ast,A}(\partial\Omega)}=\sum_{n\in\mathbb{N}}\alpha_n^2\mu_n^2.$$

By Proposition~\ref{regprop}, in particular by \eqref{bineq}, for any $f\in H_{\ast}^{1/2}(\partial\Omega)$, we have that
$$c_3^2 \sum_{n\in\mathbb{N}}\alpha_n^2\mu_n^2\leq  \|u\|^2_{H^2(\Omega)}\leq c_3^{-2} \sum_{n\in\mathbb{N}}\alpha_n^2\mu_n^2,$$
possibly for a different constant $0<c_3<1$ still depending on the a priori data only.

Let us now consider a function $v\in H^2(\Omega)$ such that $h=v|_{\partial\Omega}\in H^{1/2}_{\ast}(\partial\Omega)$.
In particular $h=\sum_{n\in\mathbb{N}}\beta_n\psi_n$ for a suitable sequence
$\{\beta_n\}_{n\in\mathbb{N}}$ of real numbers.
We call $\tilde{v}$ the solution to \eqref{Dirichlet} with boundary datum given by $h$.
Then
$$\int_{\Omega}\langle A\nabla_M u,\nabla_M v\rangle_M\, d_M=
\int_{\Omega}\langle A\nabla_M u,\nabla_M \tilde{v}\rangle_M\, d_M=
\sum_{n\in\mathbb{N}}\alpha_n\beta_n.$$
If we call, for any $n\in\mathbb{N}$, $\tilde{\beta}_n=\beta_n\mu_n$, then
$$\sup_{\sum_{n}\tilde{\beta}_n^2\leq 1} \left(\sum_{n\in\mathbb{N}}\alpha_n\frac{\tilde{\beta}_n}{\mu_n}\right)=
\sup_{\sum_{n}\tilde{\beta}_n^2\leq 1}\left(\sum_{n\in\mathbb{N}}\frac{\alpha_n}{\mu_n}\tilde{\beta}_n\right) = \left(\sum_{n\in\mathbb{N}}\frac{\alpha^2_n}{\mu^2_n}\right)^{1/2}.
$$
In other words, for any $v\in H^2(\Omega)$ with $h=v|_{\partial\Omega}\in H^{1/2}_{\ast}(\partial\Omega)$ we have
\begin{multline*}
\left|\int_{\Omega}\langle A\nabla_M u,\nabla_M v\rangle_M\, d_M\right|
\leq \|f\|_{H^{-1/2}_{\ast,A}(\partial\Omega)}\left(\sum_{n\in\mathbb{N}}\beta^2_n\mu^2_n\right)^{1/2}
\leq c_3^{-1}\|f\|_{H^{-1/2}_{\ast,A}(\partial\Omega)}\|\tilde{v}\|_{H^2(\Omega)}
\\\leq c_3^{-2}\|f\|_{H^{-1/2}_{\ast,A}(\partial\Omega)}\|h\|_{H^{3/2}(\partial\Omega)}\leq
c_4^{-1}\|f\|_{H^{-1/2}_{\ast,A}(\partial\Omega)}\|v\|_{H^2(\Omega)},
\end{multline*}
where $0<c_4<1$ is a constant still depending on the a priori data only.

Now, for any $\varphi\in L^2(\Omega)$, let $w$ be the weak solution to
\begin{equation}\label{Neumannmodified}
\left\{\begin{array}{ll}
\mathrm{div}_M(A\nabla_M w)=\varphi & \text{in }\Omega\\
\langle A\nabla_M w,\nu_M\rangle_M=c & \text{on }\partial\Omega\\
\int_{\partial\Omega}w\, d\sigma=0, &
\end{array}\right.
\end{equation}
where the constant $c$ is such that
$$\int_{\partial\Omega}c\, d\sigma=\int_{\Omega}\varphi(x)\, dx.$$
By a solution we mean $w\in H^1(\Omega)$ such that
$w|_{\partial\Omega}\in H_{\ast}^{1/2}(\partial\Omega)$ and that
$$\int_{\Omega} \langle A\nabla_M w,\nabla_M \psi\rangle_M=\int_{\partial\Omega}c\psi\, d\sigma-\int_{\Omega}\varphi(x)\psi(x)\, dx\quad\text{for any }\psi\in H^1(\Omega).$$
Still by standard regularity estimates, see for instance \cite[Chapter~2]{G}, we have that
$$\|w\|_{H^2(\Omega)}\leq C_3\|\varphi\|_{L^2(\Omega)},$$
where $C_3$ is a constant depending on the a priori data only.

We conclude that, for any $\varphi\in L^2(\Omega)$,
\begin{multline*}
\left|\int_{\Omega}u(x)\varphi(x)\, dx\right|=\left|\int_{\Omega}\langle A\nabla_M u,\nabla_M w\rangle_M\, d_M\right|\\\leq
c_4^{-1}\|f\|_{H^{-1/2}_{\ast,A}(\partial\Omega)}\|w\|_{H^2(\Omega)} \leq
C_3c_4^{-1}\|f\|_{H^{-1/2}_{\ast,A}(\partial\Omega)}\|\varphi\|_{L^2(\Omega)},
\end{multline*}
therefore
$$\|u\|_{L^2(\Omega)}\leq
C_3c_4^{-1}\|f\|_{H^{-1/2}_{\ast,A}(\partial\Omega)}$$
and the proof is concluded.\cvd

\section{The distance function from the boundary}\label{distsec}

Let $M$ be a Riemannian manifold as in Definition~\ref{manifold}.
We begin by investigating the consequences of assuming that $\varphi_M$ is smooth enough, namely we consider the following.

\begin{assum}\label{regdistassum}
For $M$, a Riemannian manifold as in Definition~\ref{manifold}, we assume that there exists $d_0>0$ such that $\varphi_M\in C^{1,1}(U_M^{d_0})$.
\end{assum}

The first consequence of Assumption~\ref{regdistassum} is the following.

\begin{prop}\label{normofdistance}
Under Assumption~\textnormal{\ref{regdistassum}}, we have
\begin{equation}\label{varphiprop}
\|\nabla_M\varphi\|_M=1\text{ in }U_M^{d_0}.
\end{equation}
\end{prop}

\proof{.} We divide the proof into several steps.

\smallskip

\noindent
\emph{First step}. We show that $\nabla\varphi_M$ is different from $0$ on $\partial\Omega$.
In fact, for any $x\in\partial\Omega$ we have
\begin{multline*}
-\frac{\partial \varphi_M}{\partial \nu(x)}(x)=\lim_{t\to 0^+}\frac{\varphi_M(x-t\nu(x))-\varphi_M(x)}{t}=
\lim_{t\to 0^+}\frac{\varphi_M(x-t\nu(x))}{t}\\ \geq \sqrt{\lambda}\lim_{t\to 0^+}\frac{\varphi(x-t\nu(x))}{t}
=
\sqrt{\lambda}\lim_{t\to 0^+}\frac{\varphi(x-t\nu(x))-\varphi(x)}{t}
=\sqrt{\lambda}>0.
\end{multline*}
In the last equality we used \eqref{varphiproper}.

\smallskip

\noindent
\emph{Second step}. We prove that $\|\nabla_M\varphi_M\|_M\leq 1$ in $U_M^{d_0}$. This follows from the obvious fact that $\varphi_M$ is Lipschitz with Lipschitz constant $1$ with respect to the distance $d_M$, that is,
$$|\varphi_M(x)-\varphi_M(y)|\leq d_M(x,y)\quad\text{for any }x,y\in \overline{\Omega}.$$

Then, 
let $x\in U_M^{d_0}$ and let $\gamma:[0,1]\to\overline{\Omega}$ be a $C^1$ curve such that $\gamma(0)=x$ and $\gamma'(0)=v$, with $v=-\nu(x)$ if $x\in\partial\Omega$. We have
$$\frac{d}{dt}(\varphi_M\circ\gamma)(0)=\nabla\varphi_M(x)v=\langle(\nabla_M\varphi_M(x))^T,v\rangle_M.$$
On the other hand,
\begin{multline*}
\left|\frac{d}{dt}(\varphi_M\circ\gamma)(0)\right|=\lim_{t\to 0^+}\frac{|\varphi_M(\gamma(t))-\varphi_M(\gamma(0))|}{t}\\
\leq \lim_{t\to 0^+}\frac{d_M(\gamma(t),\gamma(0))}{t}
\leq
\lim_{t\to 0^+}\frac{\int_0^t\|\gamma'(s)\|_M\, ds}{t}=\|\gamma'(0)\|_M=\|v\|_M.
\end{multline*}
Thus, for any $v$ or for $v=-\nu(x)$ if $x\in\partial\Omega$, we have
$$\left|\langle(\nabla_M\varphi_M(x))^T,v\rangle_M\right|\leq \|v\|_M,$$
hence $\|\nabla_M\varphi_M(x)\|_M\leq 1$.

\smallskip

\noindent
\emph{Third step}. By the first step and continuity, there exists $d_1$, $0<d_1\leq d_0$, such that we have
$0<\|\nabla_M\varphi_M(x)\|_M\leq 1$ for any $x\in U_M^{d_1}$.

We show that $\|\nabla_M\varphi_M(x)\|_M= 1$ for any $x\in U_M^{d_1}$.
By contradiction, we assume there exist $x_0\in U_M^{d_1}$, $r>0$ and $0<c<1$ such that $B_r(x_0)\subset U_M^{d_1}$ and
$0<\|\nabla_M\varphi_M(y)\|_M\leq c$ for any $y\in B_r(x_0)$.
In particular, there exists $0<t_0$ such that
$y_0=x_0+t_0\nabla\varphi_M\in B_r(x_0)$ and it satisfies the following conditions
$$
\varphi_M(y_0)>\varphi_M(x_0)\quad\text{and}\quad 2d_M(x_0,y_0)\leq d_M(y_0,y)\text{ for any }y\in U_M^{d_1}\backslash B_r(x_0).
$$
We call $h=\varphi_M(y_0)-\varphi_M(x_0)$ and we obviously have $0<h\leq d_M(x_0,y_0)$. Finally, we fix $\varepsilon$ such that
$$0<\varepsilon< \min\left(1,\frac{1-c}{c}\right) h\quad\text{and}\quad \varphi_M(y_0)+\varepsilon<d_1.$$
Let $\gamma\in \Gamma$ be such that $\gamma([0,1])\subset U^{d_1}_M$, $\gamma(0)=y_0$, $\gamma(1)\in\partial\Omega$ and
$$\mathrm{length}_M(\gamma)\leq \varphi_M(y_0)+\varepsilon<d_1.$$
There must be $s_0$, $0<s_0\leq 1$, such that $\varphi_M(\gamma(s_0))=\varphi_M(x_0)$. Therefore
$$h=\varphi_M(\gamma(0))-\varphi_M(\gamma(s_0))\leq d_M(\gamma(0),\gamma(s_0))\leq
\mathrm{length}_M(\gamma([0,s_0]))\leq h+\varepsilon.$$
But $\gamma([0,s_0]))\subset B_r(x_0)$, otherwise
$$0<2h\leq2 d_M(x_0,y_0)\leq \mathrm{length}_M(\gamma([0,s_0]))\leq h+\varepsilon<2h$$
which leads to a contradiction. Therefore,
\begin{multline*}
h=\varphi_M(\gamma(0))-\varphi_M(\gamma(s_0))=\left|-\int_0^{s_0}\nabla\varphi_M(\gamma(t))\gamma'(t)\, dt\right|\\=
\left|\int_0^{s_0}\langle(\nabla_M\varphi_M(\gamma(t)))^T,\gamma'(t)\rangle_M\, dt\right|
\leq \int_0^{s_0}c\|\gamma'(t)\|_M\, dt\\
= c\, \mathrm{length}_M(\gamma([0,s_0]))\leq c(h+\varepsilon)<h
\end{multline*}
which leads to a contradiction, thus $\|\nabla_M\varphi_M(x)\|_M= 1$ for any $x\in U_M^{d_1}$.

\smallskip

\noindent
\emph{Fourth step}. Let
$$d_2=\sup\{d:\ 0<d\leq d_0\text{ and }\|\nabla_M\varphi_M(x)\|_M= 1\text{ for any }x\in U_M^d\}.$$
By the third step we have $d_1\leq d_2$. If $d_2=d_0$ then the result is proved. Assume, by contradiction, that $d_2<d_0$. Then, by continuity, there exists $d$, $d_2<d<d_0$, such that $0<\|\nabla_M\varphi_M(x)\|_M\leq 1$ for any $x\in U_M^d$. By the same reasoning used in the third step, we conclude that 
$\|\nabla_M\varphi_M\|_M=1$ in $U_M^d$, which contradicts the definition of $d_2$.\cvd

\smallskip

Under Assumption~\ref{regdistassum}, we have that, for any $0\leq d<d_0$, 
 $\Omega_M^d$ is a $C^{1,1}$ open set and
$\partial(\Omega_M^d)=\partial\Omega_M^d$. Let $\nu$ denote the exterior normal to $\Omega_M^d$ on
$\partial\Omega_M^d$, and $\nu_M$ its corresponding one in the Riemannian setting. Then we have.

\begin{prop}\label{exteriornorm}
Under Assumption~\textnormal{\ref{regdistassum}}, for any $0\leq d<d_0$, we have
\begin{equation}\label{varphipropRiem}
(\nabla_M\varphi_M)^T=-\nu_M\text{ on }\partial\Omega_M^d.
\end{equation}
In particular this is true on $\partial\Omega$.
\end{prop}

\proof{.} It is clear that, for any $x\in\partial\Omega_M^d$, we have $(\nabla\varphi_M(x))^T=-a(x)\nu(x)$ for some positive constant
$a(x)$ depending on $x$. By the definitions of $\nabla_M\varphi_M(x)$ and of $\nu_M(x)$, we easily conclude that
$(\nabla_M\varphi_M(x))^T=-a_1(x)\nu_M(x)$ for some positive constant
$a_1(x)$ depending on $x$. Since, by Proposition~\ref{normofdistance}, $\|\nabla_M\varphi_M(x)\|_M=\|\nu_M(x)\|_M=1$,
the result immediately follows.\cvd

\smallskip

\begin{oss}\label{regUd}
Under Assumption~\ref{regdistassum}, if $\|\varphi_M\|_{C^{1,1}(U_M^{d_0})}\leq C_0$, then
$U_M^{d_0}$ is a $C^{1,1}$ open set with constants $r_1$ and $L_1$ depending on 
$r$, $L$, $R$, $d_0$ and $C_0$ only. This result can be obtained by an approximation argument, namely by
suitably approximating $\partial U_M^{d_0}\cap \Omega$ with $\partial\Omega^d_M$ as $d\to d_0^-$.
\end{oss}

The key point is the following complementary result.

\begin{prop}\label{viceversaprop}
Fixed $d_0>0$, let $f\in C^{1,1}(U_M^{d_0})$ be a nonnegative function such that
$$\|\nabla_Mf\|_M=1\text{ in }U_M^{d_0}\quad\text{and}\quad f=0\text{ on }\partial\Omega.$$
Then $f=\varphi_M$ on $U_M^{d_0}$.

Moreover, if $\|f\|_{C^{1,1}(U_M^{d_0})}\leq C_0$, we have that
\begin{equation}\label{Lipschitznabla}
\|\nabla_Mf\|_{C^{0,1}(U_M^{d_0},\mathbb{R}^N)}\leq C_1,
\end{equation}
with $C_1$ depending on $C_0$, $\lambda$ and the Lipschitz constant of the metric $G$ only.
\end{prop}

\proof{.}
First of all, we note that, since $f=0$ on $\partial\Omega$ and $f\geq 0$ in $U_M^{d_0}$, for any $x\in\partial\Omega$ we have $(\nabla f(x))^T=-a(x)\nu(x)$ for some positive constant $a(x)$ depending on $x$, thus, reasoning as in Proposition~\ref{exteriornorm}, $(\nabla_M f)^T=-\nu_M$ on $\partial\Omega$. We can also easily conclude that $f>0$ on $U_M^{d_0}\backslash\partial\Omega$.

Let $x\in U_M^{d_0}\backslash\partial\Omega$. Fixed $y\in\partial\Omega$, let $\gamma\in \Gamma(y,x)$.
Without loss of generality we can assume that $\gamma([0,1])\subset U^{d_0}_M$. Then
\begin{multline*}
f(x)=f(x)-f(y)=\int_0^1\nabla f(\gamma(s))\gamma'(s)\, ds\\ =
\int_0^1
\langle(\nabla_Mf(\gamma(s)))^T,\gamma'(s)\rangle_M\, ds\leq
\int_0^1\|\gamma'(s)\|_M\, ds=
\mathrm{length}_M(\gamma).
\end{multline*}
We can conclude that
\begin{equation}\label{ina}
f(x)\leq \varphi_M(x)\quad\text{for any }x\in U_M^{d_0}.
\end{equation}

Since $\nabla f$ is Lipschitz, by the definition of $\nabla_Mf$ and the properties of $G$, we immediately infer that also
$\nabla_M f$ is Lipschitz. Analogously, one can prove 
\eqref{Lipschitznabla}.

For any $x\in U_M^{d_0}$, let $\gamma_x$ be the (maximal) solution to the Cauchy problem for the ordinary differential equation
$$\left\{\begin{array}{l}
\gamma_x'=(\nabla_M f(\gamma_x))^T,
\\
\gamma_x(0)=x.
\end{array}\right.$$
Since $\nabla_M f$ is Lipschitz, we have existence and uniqueness of a solution $\gamma_x$ for $t\in [0,T)$, for some suitable $T>0$ depending on $x$,  even if $x\in\partial\Omega$. Moreover, for any $x\in U_M^{d_0}$, $\gamma_x(t)\in U_M^{d_0}\backslash\partial\Omega$
for any $0<t<T$.

Let $x\in U_M^{d_0}$. For any $t_0$, $t_1\in\mathbb{R}$ such that $t_0<t_1$ and for which $\gamma_x$ is defined, let us call $z_0=\gamma_x(t_0)$ and $z_1=\gamma_x(t_1)$. Then we observe that
\begin{multline*}
f(z_1)-f(z_0)=\int_{t_0}^{t_1}\nabla f(\gamma_x(s))\gamma_x'(s)\, ds\\ =
\int_{t_0}^{t_1}
\langle(\nabla_Mf(\gamma_x(s)))^T,\gamma_x'(s)\rangle_M\, ds
=
\mathrm{length}_M(\gamma_x([t_0,t_1]))=t_1-t_0,
\end{multline*}
therefore $d_M(z_0,z_1)\leq f(z_1)-f(z_0)$. In particular, if $z_0\in\partial\Omega$, then
$$\varphi_M(z_1)\leq d_M(z_0,z_1)\leq f(z_1)-f(z_0)=f(z_1),$$
thus, by the previous inequality \eqref{ina}, we have $\varphi_M(z_1)=f(z_1)$.

We claim the following result. Let $d_1$, $0\leq d_1<d_0$ be such that $f(x)=\varphi_M(x)$ for any $x\in U_M^{d_0}$ with $\varphi_M(x)\leq d_1$. Then there exists $d$  such that
$d_1<d< d_0$ and $f(x)=\varphi_M(x)$ for any $x\in U_M^{d_0}$ with $\varphi_M(x)\leq d$.

In order to prove the claim, let us begin with the following remark, where we assume that $d_1>0$. Let $x\in U_M^{d_0}$ be such that $\varphi_M(x)=f(x)=d_1$. By the implicit function theorem, there exist a $C^1$
function $\phi_x:\mathbb{R}^{N-1}\to\mathbb{R}$ and an open neighbourhood $U_x$ of $x$
such that for any $y\in U_x$ we have, up to a rigid transformation depending on $x$,
\begin{equation}\label{levelf}
\left\{y\in U_x:\ f(y)\lesseqqgtr d_1\right\}=\left\{y=(y',y_N)\in U_x:\  y_N\lesseqqgtr\phi_x(y')\right\}.
\end{equation}
Without loss of generality, up to changing $U_x$, we can assume that
$$U_x^-=\left\{y=(y',y_N)\in U_x:\  y_N<\phi_x(y')\right\}$$ is connected.
We want to show that \eqref{levelf} holds true even if we replace $f$ with $\varphi_M$. By \eqref{ina},
it is clear that
$$\left\{y\in U_x:\ \varphi_M(y)> d_1\right\}\supset\left\{y\in U_x:\ f(y)> d_1\right\}=\left\{y=(y',y_N)\in U_x:\  y_N>\phi_x(y')\right\}.$$
Moreover, by our assumption,
$$\left\{y\in U_x:\ \varphi_M(y)\leqq d_1\right\}\subset\left\{y\in U_x:\ f(y)\leqq d_1\right\}=\left\{y=(y',y_N)\in U_x:\  y_N\leqq\phi_x(y')\right\}.$$
Since $\varphi_M$ can not have interior local minimum points, there exists $y_1\in U_x$ such that $\varphi_M(y_1)< d_1$.
Then $f(y_1)<d_1$ and $y_1\in U_x^-$.
Assume by contradiction that there exists $y_2\in U_x$ such that
$\varphi_M(y_2)>d_1\geq f(y_2)$. Actually, by continuity, we can always assume that 
$\varphi_M(y_2)>d_1> f(y_2)$, hence that $y_2\in U^-_x$ as well. We connect $y_1$ to $y_2$ with a smooth curve all contained in $U_x^-$. There must be a point $y$ along this curve on which $\varphi_M(y)=d_1$, thus we obtain a contradiction since $f(y)<d_1$.

This remark allows us to show that there exists $\varepsilon>0$ such that $d_1+\varepsilon<d_0$ and
$f(y)>d_1$ for any $y$ with $d_1<\varphi_M(y)<d_1+\varepsilon$. Assume by contradiction that there exists $y_0$ such that
$d_1<f(y_0)<\varphi_M(y_0)\leq  d_1+\varepsilon/4$. We note that it is well-defined $z_0=\gamma_{y_0}(d_1-f(y_0))$. It happens that $f(z_0)=d_1$ and $d_M(z_0,y_0)\leq (f(y_0)-d_1)\leq \varepsilon/4$, therefore $\varphi_M(z_0)\leq d_1+\varepsilon/2$. By \eqref{ina}, $d_1\leq \varphi_M(z_0)$ but
$\varphi_M(z_0)$ can not be greater than $d_1$, otherwise $f(z_0)$ should be greater than $d_1$ as well.
We conclude that $\varphi_M(z_0)=d_1$, therefore
$\varphi_M(y_0)\leq \varphi_M(z_0)+d_M(z_0,y_0)\leq d_1+(f(y_0)-d_1)=f(y_0)$ which gives the contradiction and proves the claim.

Let us conclude the proof by defining
$$d_2=\sup\{d:\ 0<d< d_0\text{ and }f(x)=\varphi_M(x)\text{ for any }x\in U_M^{d_0}\text{ with }\varphi_M(x)\leq d\}.$$
If $d_2=d_0$ the proof is concluded. If, by contradiction, $d_2<d_0$, by continuity we have that $f(x)=\varphi_M(x)$ for any $x\in U_M^{d_0}$ with $\varphi_M(x)\leq d_2$ and the claim contradicts the definition of $d_2$.\cvd

\smallskip

We point out the following important property. Under Assumption~\ref{regdistassum}, or equivalently under the assumptions of Proposition~\ref{viceversaprop}, for any $x\in U_M^{d_0}$, let $\gamma_x$ be the (maximal) solution to the Cauchy problem for the ordinary differential equation
\begin{equation}\label{ode}
\left\{\begin{array}{l}
\gamma_x'=(\nabla_M \varphi_M(\gamma_x))^T,
\\
\gamma_x(0)=x.
\end{array}\right.
\end{equation}
Then $\gamma_x:[-\varphi_M(x),d_0-\varphi_M(x))$
with $\gamma_x(-\varphi_M(x))=y\in \partial\Omega$. In other words, for any $x\in U_M^{d_0}$ there exists $y\in\partial\Omega$
such that $x=\gamma_y(\varphi_M(x))$ and
$$\varphi_M(x)=d_M(x,y)=\mathrm{length}_M(\gamma_y([0,\varphi_M(x)])).$$
We can then state the following result.

\begin{cor}\label{localcoordinate}
Under Assumption~\textnormal{\ref{regdistassum}}, or equivalently under the assumptions of Proposition~\textnormal{\ref{viceversaprop}},
we can define a coordinate system for $U_M^{d_0}$ given by $T:\partial\Omega\times [0,d_0)\to U_M^{d_0}$ such that for any $(y,d)\in \partial\Omega\times [0,d_0)$ we have $T(y,d)=\gamma_y(d)$. We note that, for any $0\leq d<d_0$, we have
$T(\partial\Omega\times{d})=\partial\Omega^d_M$.

Moreover, if we assume that $\|\varphi_M\|_{C^{1,1}(U_M^{d_0})}\leq C_0$, then
$T$ is bi-Lipschitz, that is, $T$ and its inverse $T^{-1}$ are Lipschitz, with Lipschitz constants bounded by a constant depending on $C_0$, $d_0$, $\lambda$, the Lipschitz constant of the metric $G$ and $C(\Omega)$ as in \eqref{lipd} only.
\end{cor}

\proof{.}
The fact that $T$ is injective simply depends on the uniqueness for the solution to \eqref{ode}. We begin by showing that $T$ is Lipschitz, using an argument that is related to the continuity of solutions to ordinary differential equations with respect to the data.

First of all, as for \eqref{Lipschitznabla}, we note that
\begin{equation}\label{Lipschitznabla2}
\|\nabla_M\varphi_M\|_{C^{0,1}(U_M^{d_0},\mathbb{R}^N)}\leq C_1,
\end{equation}
with $C_1$ depending on $C_0$, $\lambda$ and the Lipschitz constant of the metric $G$ only.

For any $i=1,2$, let 
$x_i\in U_M^{d_0}$ and $t_i\in [-\varphi_M(x_i),d_0-\varphi_M(x_i))$. We wish to estimate
$\|\gamma_{x_2}(t_2)-\gamma_{x_1}(t_1)\|$. By Volterra integral equation, we have that
$$\gamma_{x_2}(t_2)-\gamma_{x_1}(t_1)=
\left(x_2+\int_0^{t_2}\nabla_M \varphi_M(\gamma_{x_2}(s))\, ds\right) -\left(x_1+\int_0^{t_1}\nabla_M \varphi_M(\gamma_{x_1}(s))\, ds\right).
$$

We begin by considering the case $t_1=t_2$.
Then \begin{multline}\label{uno}
\|\gamma_{x_2}(t_1)-\gamma_{x_1}(t_1)\|\leq\|x_2-x_1\|+
\left|\int_0^{t_1}\| \nabla_M \varphi_M(\gamma_{x_2}(s))-\nabla_M \varphi_M(\gamma_{x_1}(s)) \|\, ds\right|
\\\leq
\|x_2-x_1\|+C_1\left|\int_0^{t_1}\| \gamma_{x_2}(s)-\gamma_{x_1}(s) \|\, ds\right|,
\end{multline}
where we used \eqref{Lipschitznabla2}. Then, by Gronwall lemma, we have that
\begin{equation}\label{due}
\|\gamma_{x_2}(t_1)-\gamma_{x_1}(t_1)\|\leq e^{C_1d_0}\|x_2-x_1\|.
\end{equation}
Moreover, we infer that
\begin{equation}\label{tre}
\|(\gamma_{x_2}(t_1)-x_2)-(\gamma_{x_1}(t_1)-x_1)\|\leq C_1e^{C_1d_0}\|x_2-x_1\||t_1|,
\end{equation}
an inequality that will be crucial later on.

We now turn to the general case. If $t_1\leq 0\leq t_2$, or $t_2\leq 0\leq t_1$, then
\begin{multline}\label{Lipgammaa}
\|\gamma_{x_2}(t_2)-\gamma_{x_1}(t_1)\|\leq
\|\gamma_{x_2}(t_2)-x_2\|+\|x_2-x_1\|+\|x_1-\gamma_{x_1}(t_1)\|\\\leq
\|x_2-x_1\|+\left|\int_0^{t_2}\|\nabla_M \varphi_M(\gamma_{x_2}(s))\|\, ds\right|+
\left|\int_0^{t_1}\|\nabla_M \varphi_M(\gamma_{x_1}(s))\|\, ds\right|\\\leq
\|x_2-x_1\|+\sqrt{\lambda^{-1}}(|t_1|+|t_2|)=\|x_2-x_1\|+\sqrt{\lambda^{-1}}|t_2-t_1|,
\end{multline}
where we used \eqref{normcomparison} and the fact that $\|\nabla_M \varphi_M\|_M=1$.

Otherwise, up to swapping $x_1$ with $x_2$, we have $0\leq t_1\leq t_2$ or $t_2\leq t_1\leq 0$, and then
\begin{multline}\label{Lipgamma0}
\|\gamma_{x_2}(t_2)-\gamma_{x_1}(t_1)\|\leq
\|\gamma_{x_2}(t_2)-\gamma_{x_2}(t_1)\|
+\|\gamma_{x_2}(t_1)-\gamma_{x_1}(t_1)\|\\\leq
\left|\int_{t_1}^{t_2}\|\nabla_M \varphi_M(\gamma_{x_2}(s))\|\, ds\right|+
\|\gamma_{x_2}(t_1)-\gamma_{x_1}(t_1)\|\\\leq 
\sqrt{\lambda^{-1}} |t_2-t_1|+\|\gamma_{x_2}(t_1)-\gamma_{x_1}(t_1)\|.
\end{multline}
By \eqref{due} and \eqref{Lipgamma0} we can conclude that
\begin{equation}\label{Lipgamma}
\|\gamma_{x_2}(t_2)-\gamma_{x_1}(t_1)\|\leq
e^{C_1d_0}\|x_2-x_1\|+\sqrt{\lambda^{-1}} |t_2-t_1|.
\end{equation}

By \eqref{Lipgammaa} and \eqref{Lipgamma}, it is immediate to prove that $T$ is Lipschitz and that its Lipschitz constant is bounded by a constant depending on $C_0$, $d_0$, $\lambda$ and the Lipschitz constant of the metric $G$ only.

Let us now pass to the properties of $T^{-1}$. For any $x\in U_M^{d_0}$, we have that
$$T^{-1}(x)=(\gamma_{x}(-\varphi_M(x)),\varphi_M(x))\in\partial\Omega\times [0, d_0).$$
We recall that
$\varphi_M$ is Lipschitz, with Lipschitz constant $1$, with respect to the distance $d_M$. Hence we can conclude the proof using again \eqref{Lipgamma}.\cvd

\smallskip

The following technical proposition is a crucial ingredient for the proof of our main decay estimate and it may be of independent interest as well.

\begin{prop}\label{normalderivativelemma}
Under Assumption~\textnormal{\ref{regdistassum}}, or equivalently under the assumptions of Proposition~\textnormal{\ref{viceversaprop}}, let $\|\varphi_M\|_{C^{1,1}(U_M^{d_0})}\leq C_0$.

Let $w\in W^{1,1}_{\mathrm{loc}}(\Omega)$ and let, for any $0<d<d_0$,
$$S(d)=\int_{\partial\Omega^d_M}w(x)\, d\sigma_M(x).$$

We have that $S$ is absolutely continuous on any compact subinterval of $(0,d_0)$ and, for almost any $d$, $0<d<d_0$,
\begin{multline}\label{derivative}
S'(d)=-\int_{\partial\Omega^d_M}\nabla w(x)\nu_M(x)\, d\sigma_M(x)+A(d)\\=
-\int_{\partial\Omega^d_M}\langle \nabla_M w(x),\nu_M(x)\rangle_M\, d\sigma_M(x)+A(d)
\end{multline}
where
$$|A(d)|\leq C\int_{\partial\Omega^d_M}|w(x)|\, d\sigma_M(x)$$
for a constant $C$ depending on $C_0$, $d_0$, $\lambda$ and the Lipschitz constant of the metric $G$ only.

In particular, if $w\geq 0$, then
\begin{equation}\label{derivativeresto}
|A(d)|\leq C S(d).
\end{equation}
\end{prop}

\begin{oss}\label{W11remark}
If $w\in W^{1,1}(\Omega)$, then we can define
$$S(0)=\int_{\partial\Omega^0_M}w(x)\, d\sigma_M(x)=\int_{\partial\Omega}w(x)\, d\sigma_M(x),$$
and we have that $S$ is absolutely continuous on any compact subinterval of $[0,d_0)$.
\end{oss}

\proof{.}
We just assume $w\in W^{1,1}(\Omega)$ as in Remark~\ref{W11remark},
since, when $w\in W^{1,1}_{\mathrm{loc}}(\Omega)$, the result 
easily follows by the arguments we present in the sequel.

We begin by observing that, for any $s$, $0\leq s<d_0$, we have
$$\int_{\partial\Omega^s_M}w(x)\, d\sigma_M(x)=
\int_{\partial\Omega^s_M}w(x)h(x)\, d\sigma(x)
$$
where
$$h(x)=\sqrt{\langle G^{-1}(x)\nu(x),\nu(x)\rangle}\sqrt{g(x)}.$$

Moreover, for any $s_1$, $s_2\in [0,d_0)$, we call $T_{s_1,s_2}:\partial\Omega^{s_1}_M\to\partial\Omega^{s_2}_M$ the change of coordinates such that
$$T_{s_1,s_2}(x)=\gamma_x(s_2-s_1).$$
By \eqref{due} and the fact that $T_{s_1,s_2}$ is invertible with $T^{-1}_{s_1,s_2}=T_{s_2,s_1}$, we deduce that $T_{s_1,s_2}$ is bi-Lipschitz, therefore
$$\int_{\partial\Omega^{s_2}_M}w(z)\, d\sigma_M(z)=
\int_{\partial\Omega^{s_1}_M}w(\gamma_x(s_2-s_1))h(\gamma_x(s_2-s_1))k(x)\, d\sigma(x)
$$
where $k(x)$ can be computed as follows. For almost every $x\in \partial\Omega^{s_1}_M$, with respect to the $(N-1)$-dimensional Hausdorff measure, $T_{s_1,s_2}$ admits a tangential differential at $x$. Namely, for any orthonormal basis $v_1,\ldots,v_{N-1}$ of the tangent space to $\partial\Omega^{s_1}_M$ at $x$, there exists
$$J_{\tau}(x)=J_{\tau}T_{s_1,s_2}(x)=\left[\frac{\partial T_{s_1,s_2}}{\partial v_1}(x)\cdots\frac{\partial T_{s_1,s_2}}{\partial v_{N-1}}(x)\right].$$
Then 
\begin{equation}\label{quattro}
k(x)=\sqrt{\det \left((J_{\tau}(x))^TJ_{\tau}(x)\right)}.
\end{equation}
Let us call $\tilde{T}_{s_1,s_2}=T_{s_1,s_2}-Id$ and let, analogously, $\tilde{J}_{\tau}(x)=J_{\tau}\tilde{T}_{s_1,s_2}(x)$. By \eqref{tre}, we infer that for any $i=1,\ldots,N-1$,
\begin{equation}\label{cinque}
\left\|\frac{\partial \tilde{T}_{s_1,s_2}}{\partial v_i}(x) \right\|\leq C_1e^{C_1d_0}|s_2-s_1|.
\end{equation}

Therefore, for almost every $x\in \partial\Omega^{s_1}_M$, again with respect to the $(N-1)$-dimensional Hausdorff measure, we call $a(x,s_1,s_2)$ the number such that
$$\frac{h(\gamma_x(s_2-s_1))}{h(x)}k(x)=1+a(x,s_1,s_2).$$

By using \eqref{quattro} and \eqref{cinque} to handle $k(x)$, it is not difficult to show that, for some constant $C_2$ depending on 
$C_0$, $d_0$, $\lambda$ and the Lipschitz constant of the metric $G$ only,
\begin{equation}\label{sei}
|a(x,s_1,s_2)|\leq C_2|s_2-s_1|\quad\text{for almost every }x\in \partial\Omega^{s_1}_M.
\end{equation}

Then, for almost every $x\in \partial\Omega^{s_1}_M$, or for almost every $z=\gamma_x(s_2-s_1)\in \partial\Omega^{s_2}_M$,
\begin{multline*}
w(\gamma_x(s_2-s_1))=w(x)+\int_{0}^{s_2-s_1}\nabla w(\gamma_x(s))\gamma_x'(s)\, ds\\=
w(x)-\int_{0}^{s_2-s_1}\nabla w(\gamma_x(s))\nu_M(\gamma_x(s))\, ds=
w(x)-\int_{0}^{s_2-s_1}\nabla w(\gamma_z(-s))\nu_M(\gamma_z(-s))\, ds.
\end{multline*}
We call $\Omega_{s_1,s_2}$ the following set
$$\Omega_{s_1,s_2}=\left\{\begin{array}{ll}
\Omega^{s_1}_M\backslash\Omega^{s_2}_M &\text{if }s_1\leq s_2\\
\\
\Omega^{s_2}_M\backslash\Omega^{s_1}_M &\text{if }s_2\leq s_1
\end{array}
\right.$$
and we call
$$
b(s_1,s_2)=\left\{\begin{array}{ll}
1 &\text{if }s_1< s_2\\
0 &\text{if }s_1= s_2\\
-1 &\text{if }s_1> s_2.\end{array}
\right.$$

Then, by Fubini theorem and the coarea formula,
\begin{multline*}
\int_{\partial\Omega^{s_2}_M}w(x)\, d\sigma_M(x)
-\int_{\partial\Omega^{s_1}_M}w(x)\, d\sigma_M(x)\\=
\int_{\partial\Omega^{s_1}_M}w(x)a(x,s_1,s_2)\, d\sigma_M(x)-
\int_{\partial\Omega^{s_2}_M}\left(\int_{0}^{s_2-s_1}\nabla w(\gamma_z(-s))\nu_M(\gamma_z(-s))\, ds\right)\, d\sigma_M(z)
\\=
\int_{\partial\Omega^{s_1}_M}w(x)a(x,s_1,s_2)\, d\sigma_M(x)-
\int_{s_1}^{s_2}\left(\int_{\partial\Omega^t_M}\nabla w(x)\nu_M(x)(1+a(x,t,s_2))\, d\sigma_M(x)\right)\, dt
\\=
\int_{\partial\Omega^{s_1}_M}w(x)a(x,s_1,s_2)\, d\sigma_M(x)-b(s_1,s_2)
\int_{\Omega_{s_1,s_2}}\nabla w(x)\nu_M(x)(1+a(x,s,s_2))\, d_M(x)\\=
\int_{\partial\Omega^{s_1}_M}w(x)a(x,s_1,s_2)\, d\sigma_M(x)\\-b(s_1,s_2)
\int_{\Omega_{s_1,s_2}}\nabla w(x)\nu_M(x)\, d_M(x)-b(s_1,s_2)
\int_{\Omega_{s_1,s_2}}\nabla w(x)\nu_M(x)a(x,s,s_2)\, d_M(x)
\\=A(s_1,s_2)-B(s_1,s_2)-C(s_1,s_2)
\end{multline*}
where, for any $x\in \Omega_{s_1,s_2}$ we set $s=\varphi_M(x)$.
First of all, we deduce that
$$[0,d_0)\ni s\mapsto \int_{\partial\Omega^{s}_M}w(x)\, d\sigma_M(x)$$
is a continuous function.

Again by coarea formula, we have that the function
\begin{multline*}
[0,d_0)\ni s\mapsto B(d_0/2,s)=b(d_0/2,s)\int_{\Omega_{d_0,s}}\nabla w(x)\nu_M(x)\, d_M(x)\\=
\int_{d_0/2}^{s}\left(\int_{\partial\Omega^t_M}\nabla w(x)\nu_M(x)\, d\sigma_M(x)\right)\, dt
\end{multline*}
is absolutely continuous, with respect to $s$, on any compact subinterval of $[0,d_0)$ and,
for almost every $s_1\in (0,d_0)$, we have
\begin{multline*}
B'(d_0/2,s_1)=\lim_{s_2\to s_1}\frac{B(d_0/2,s_2)-B(d_0/2,s_1)}{s_2-s_1}\\=
\lim_{s_2\to s_1}\frac{B(s_1,s_2)}{s_2-s_1}
=\int_{\partial\Omega^{s_1}_M}\nabla w(x)\nu_M(x)\, d\sigma_M(x).
\end{multline*}

The function
$$[0,d_0)\ni s\mapsto D(s)=\int_{\partial\Omega^{s}_M}w(x)\, d\sigma_M(x)+B(d_0/2,s)
$$
is clearly Lipschitz continuous on any compact subinterval of $[0,d_0)$, therefore, for almost every $s_1\in (0,d_0)$, there exists
$$
D'(s_1)=\lim_{s_2\to s_1}\frac{D(s_2)-D(s_1)}{s_2-s_1}=
\lim_{s_2\to s_1}\frac{A(s_1,s_2)-C(s_1,s_2)}{s_2-s_1}.
$$
It is easy to see that
$$\frac{C(s_1,s_2)}{s_2-s_1}\to 0\quad\text{as }s_2\to s_1$$
and that
$$\left|\frac{A(s_1,s_2)}{s_2-s_1}\right|\leq C_2\int_{\partial\Omega^{s_1}_M}|w(x)|\, d\sigma_M(x).$$
Therefore the proof can be easily concluded.\cvd

\smallskip

Our aim is to modify our metric $G$ near the boundary of $\Omega$, by multiplying it with a scalar function $\eta$, in such a way that the new metric satisfies Assumption~\ref{regdistassum}. 
The construction is given in the next theorem.

\begin{teo}\label{AKSmethod}
Let us fix positive constants $R$, $r$ and $L$.
Let $\Omega\subset\overline{B_R(0)}\subset\mathbb{R}^N$ be a bounded open set of class $C^{1,1}$ with constants $r$ and $L$. Let us consider $\tilde{d}_0>0$ as in Theorem~\textnormal{\ref{Del-Zolthm}} and $\varphi$ the distance to the boundary of $\Omega$ as in Definition~\textnormal{\ref{Eucldistfun}}.

Let  $G$ be a Lipschitz symmetric tensor in $\Omega$ which is uniformly elliptic with constant $\lambda$, $0<\lambda<1$, in $\Omega$
and such that $\|G\|_{C^{0,1}(\overline{\Omega})}\leq C$.

Then there exist a constant $C_1>0$, depending on $r$, $L$, $R$, $\lambda$ and $C$ only, and a function
 $\eta\in C^{0,1}(\overline{\Omega})$, which is uniformly elliptic with constant $\lambda$ in $\Omega$ and such that
 $\|\eta\|_{C^{0,1}(\overline{\Omega})}\leq C_1$, such that the following holds.
 
 Let us call $\tilde{G}=\eta G$ and $\tilde{M}$ the corresponding Riemannian manifold on $\overline{\Omega}$. Let $\varphi_{\tilde{M}}$ be the corresponding distance from the boundary and, for any $d\geq 0$, $U_{\tilde{M}}^d=\{x\in\overline{\Omega}:\ \varphi_{\tilde{M}}(x)<d\}$.
 
 Then we have that $U^{\tilde{d}_0/2}=U_{\tilde{M}}^{\tilde{d}_0/2}$ and
\begin{equation}\label{crucialequality}
\varphi_{\tilde{M}}=\varphi\quad\text{in }U_{\tilde{M}}^{\tilde{d}_0/2}.
\end{equation}
 \end{teo}

\proof{.} Let us define $\hat{\eta}:U^{\tilde{d}_0}\to\mathbb{R}$ such that
$$\hat{\eta}=\|\nabla_M\varphi\|_M^2\text{ in }U^{\tilde{d}_0}.$$
By \eqref{normcomparison}, we obtain that $\lambda\leq \hat{\eta}\leq\lambda^{-1}$ in $U^{\tilde{d}_0}$, and we have that 
$$\hat{\eta}^{-1}\|\nabla_M\varphi\|_M^2=1\text{ in }U^{\tilde{d}_0}.$$

Then we fix a cutoff function $\chi\in C^{\infty}(\mathbb{R})$ such that $\chi$ is decreasing, $\chi(t)=1$ for any $t\leq \tilde{d}_0/2$ and $\chi(t)=0$ for any $t\geq 3\tilde{d}_0/4$. We define, for any $x\in\overline{\Omega}$,
$$\eta(x)=
\chi(\varphi(x))\hat{\eta}(x)+(1-\chi(\varphi(x))$$
and we observe that $\lambda\leq \eta\leq\lambda^{-1}$ in $\overline{\Omega}$.

Let $\tilde{G}=\eta G$. By construction of $\eta$ and by \eqref{gradcomp},
we have that
$$\|\nabla_{\tilde{M}}\varphi\|_{\tilde{M}}=1\quad\text{in }U^{\tilde{d}_0/2}.$$
Therefore, applying Proposition~\ref{viceversaprop} with $f=\varphi$, we conclude that, at least in a neighbourhood
of $\partial\Omega$, $\varphi_{\tilde{M}}=\varphi$. It is not difficult to show that
such a neighbourhood is actually equal to $U^{\tilde{d}_0/2}$ and that it coincides with $U_{\tilde{M}}^{\tilde{d}_0/2}$
as well.

It remains to show the Lipschitz regularity of $\eta$ and for this purpose it is enough to show that $\hat{\eta}$ is Lipschitz
in $U^{\tilde{d}_0}$. Again by \eqref{gradcomp}, we infer that for any $x\in U^{\tilde{d}_0}$
$$\hat{\eta}(x)=\|\nabla_M\varphi(x)\|_M^2=\langle (\nabla \varphi(x))^T,G^{-1}(x)(\nabla \varphi(x))^T\rangle.$$
Then we can easily conclude by exploiting the Lipschitz regularity of $G$ and the fact that $\varphi\in C^{1,1}(U^{\tilde{d}_0})$ as proved in Therem~\ref{Del-Zolthm}.\cvd

\smallskip

We conclude that $\tilde{G}=\eta G$ constructed in Theorem~\ref{AKSmethod} is a Lipschitz symmetric tensor in $\Omega$ which is uniformly elliptic with constant $\lambda_1=\lambda^2$ in $\Omega$ and
such that $\|\tilde{G}\|_{C^{0,1}(\overline{\Omega})}\leq C_2$, with $C_2$ depending on $C$, $C_1$ and $\lambda$ only. Moreover, by Theorem~\ref{Del-Zolthm} and \eqref{crucialequality}, $\tilde{G}$ satisfies Assumption~\ref{regdistassum} with $d_0=\tilde{d}_0/2$.

 \section{The decay estimate}\label{decaysec}
 
Let us fix positive constants $R$, $r$, $L$, $C_0$, $C_1$, $\lambda$ and $\lambda_1$, 
with $0<\lambda<1$ and $0<\lambda_1<1$. We refer to these constants as the
\emph{a priori data}.

Let $\Omega\subset\overline{B_R(0)}\subset\mathbb{R}^N$ be a bounded domain of class $C^{1,1}$ with constants $r$ and $L$. Let us consider $\tilde{d}_0>0$ as in Theorem~\textnormal{\ref{Del-Zolthm}} and $\varphi$ the distance to the boundary of $\Omega$ as in Definition~\textnormal{\ref{Eucldistfun}}.

Let  $G$ be a Lipschitz symmetric tensor in $\Omega$ which is uniformly elliptic with constant $\lambda$
and such that $\|G\|_{C^{0,1}(\overline{\Omega})}\leq C_0$.

Let $A$ be a Lipschtitz conductivity tensor in $\Omega$ 
which is uniformly elliptic with constant $\lambda_1$
and such that
$\|A\|_{C^{0,1}(\overline{\Omega})}\leq C_1$.

We further suppose that Assumption~\ref{scalarassum} holds.

Let us fix $f\in H^{1/2}(\partial\Omega)$, with $f\neq 0$, and let us call $\Phi$ its frequency as in Definition~\ref{frequencydefin}.
We assume that $\Phi>0$, that is, $f$ is not constant on $\partial\Omega$.
Let $u\in H^1(\Omega)$ be the solution to \eqref{Dirichlet}. We recall that $u\in H^2_{\mathrm{loc}}(\Omega)$ and the equation is satisfied pointwise almost everywhere in $\Omega$.

The important remark is that, without loss of generality, we can assume that the following fact holds.

By Remark~\ref{anisotropy}, we can assume that
\begin{equation}
\label{isoassum}
A=\gamma I_N\text{ with }\gamma\in C^{0,1}(\overline{\Omega}).
\end{equation}

We can assume that
$G$ satisfies Assumption~\ref{regdistassum} with some positive constant $d_0$. Under this assumption, we need to add $d_0$ and $\|\varphi_M\|_{C^{1,1}(U_M^{d_0})}$ to the a priori data. In particular, by Theorem~\ref{AKSmethod} and Remark~\ref{conformalchange}, 
 we can assume that 
 \begin{equation}\label{crucialequality2}
d_0=\tilde{d}_0/2,\quad U^{d_0}=U_M^{d_0}\quad\text{and}\quad
\varphi_M=\varphi\quad\text{in }U_M^{d_0}.
\end{equation}
In this case, by Theorem~\ref{Del-Zolthm}, $d_0$ and $\|\varphi_M\|_{C^{1,1}(U_M^{d_0})}$ depend on $r$, $L$ and $R$ only.

Before stating our decay estimates, we need to set some notation.
For any $0\leq d<d_0$, let us define
$$D(d)=\int_{\Omega^d_M}\gamma(x)\|\nabla_M u(x)\|_M^2\, d_M(x)\quad\text{and}\quad H(d)=\int_{\partial\Omega^d_M}\gamma(x) u^2(x)\, d\sigma_M(x).$$
We recall that, for any such $d$, $\partial\Omega^d_M=\partial(\Omega^d_M)$ and, if \eqref{crucialequality2} holds,
$\Omega^d=\Omega^d_M$ and $\partial\Omega^d=\partial\Omega^d_M=\partial(\Omega^d)$. Moreover, by unique continuation, for example by \cite{Gar-Lin1} for $N\geq 3$, and the maximum principle, both $D(d)$ and $H(d)$ must be strictly positive for any $0\leq d<d_0$. 
We define the \emph{frequency function} $N$ as follows
\begin{equation}
N(d)=\frac{D(d)}{H(d)},\qquad 0\leq d<d_0.
\end{equation}
We note that, by Remark~\ref{freqrem},
there exists a constant $c_1$, $0<c_1<1$ depending on $\lambda$ and $\lambda_1$ only, such that
\begin{equation}
\lambda_1c_1\Phi\leq N(0)\leq (\lambda_1c_1)^{-1}\Phi
\end{equation}
where $\Phi$ is the frequency of the boundary datum $f$.

For any $s\geq 0$ we define
\begin{equation}\label{hfundef}
h(s)=\left\{\begin{array}{ll}e^{-s} &\text{if }s\leq 1\\
(es)^{-1}&\text{if }s>1.
\end{array}\right.
\end{equation}
We note that $h(0)=1$ and $h$ is a positive $C^1$ strictly decreasing function.

\begin{teo}\label{mainthm}
Let $f\in H^{1/2}(\partial\Omega)$, with $f\neq 0$, and let its frequency $\Phi$ be positive. 
Under the previous assumptions and notation, there exist two positive constants $C_2$ and $c_2$, depending on the a priori data only,
such that, for any $d$, $0<d<d_0$, we have
\begin{equation}\label{decayestimate}
D(d)\leq e^{C_2d}D(0)h(c_2d\Phi).
\end{equation}
\end{teo}

In the next theorem, we control the decay of the function, instead of that of its gradient.
Namely, we assume that $f\in H^{1/2}_{\ast}(\partial\Omega)$, with $f\neq 0$, and that $\Phi_1$ is its lower frequency. We recall that $\Phi_1\leq \Phi$.

\begin{teo}\label{mainthmbis}
Let $f\in H^{1/2}_{\ast}(\partial\Omega)$, with $f\neq 0$, and let $\Phi_1$ be its lower frequency. 
Under the previous assumptions and notation, there exist two positive constants $C_3$ and $c_3$, depending on the a priori data only,
such that, for any $d$, $0<d<d_0/2$, we have
\begin{equation}\label{decayestimatebis}
H(d)\leq e^{C_3d}H(0)h(c_3d\Phi_1).
\end{equation}
\end{teo}

As a corollary, we obtain a higher order decay for $D$ with respect to the lower frequency.

\begin{cor}\label{maincor}
Let $f\in H^{1/2}_{\ast}(\partial\Omega)$, with $f\neq 0$, and let $\Phi_1$ be its lower frequency. 
Under the previous assumptions and notation, there exists a further absolute positive constant $C_4$
such that, for any $d$, $0<d<d_0/4$, we have
\begin{multline}\label{decayestimateter}
D(d)\leq \frac{C_4}{d}e^{3C_3d/2}H(0)h(c_3d\Phi_1/2)\\\leq C_4e^{3C_3d/2}D(0)\frac{h(c_3d\Phi_1/2)}{\lambda_1c_1d\Phi}\leq
C_4e^{3C_3d/2}D(0)\frac{h(c_3d\Phi_1/2)}{\lambda_1c_1d\Phi_1}.
\end{multline}
\end{cor}

\begin{oss}\label{Steklovremark}
If $f=\phi$ where $\phi$ is a Steklov eigenfunction with Steklov eigenvalue $\mu>0$, that is, $u=w$ where $w$ is a solution to \eqref{Steklov}, then the results of Theorems~\ref{mainthm} and \ref{mainthmbis} and of Corollary~\ref{maincor} still hold, possibly with different constants still depending on the a priori data only, even if we replace both $\Phi$ and $\Phi_1$ with the Steklov eigenvalue $\mu$.  
\end{oss}

\proof{ of Corollary~\textnormal{\ref{maincor}}.} We sketch the proof of the corollary. Let $0<d<d_0/4$. Then we have, by coarea formula and \eqref{decayestimatebis},
\begin{equation}\label{bbb}
\int_{\Omega_M^{d/2}\backslash\Omega_M^{3d/2}}\gamma u^2\, d_M
=\int_{d/2}^{3d/2}H(t)\, dt
\leq H(0)
de^{3C_3d/2}h(c_3d\Phi_1/2).
\end{equation}

Then we apply a Caccioppoli inequality. Let $\chi\in C^{\infty}_0(\mathbb{R})$ be an even positive function such that
$\chi$ is decreasing on $[0,1)$, $\chi=1$ on $[0,1/2]$ and $\chi=0$ on $[3/4,+\infty)$. We define the function $\eta_d$
as follows
$$\eta_d(x)=\chi\left(2\frac{\varphi_M(x)-d}{d}\right)\quad\text{for any }x\in\Omega$$
and we note that
$$\nabla_M\eta_d(x)=\frac{2}{d}\chi'\left(2\frac{\varphi_M(x)-d}{d}\right)\nabla_M\varphi_M(x)
\quad\text{for any }x\in\Omega.$$
Therefore
$$\|\nabla_M\eta_d(x)\|_M\leq \frac{C}{d}\quad\text{for any }x\in\Omega\backslash\Omega_M^{d_0/2}$$
where $C$ is an absolute constant.

Then
\begin{multline*}
0=\int_{\Omega}\gamma\langle \nabla_Mu,\nabla_M(u\eta_d^2)\rangle_M\, d_M\\=
\int_{\Omega}\gamma\langle \nabla_Mu,\nabla_Mu\rangle_M\eta_d^2\, d_M+2
\int_{\Omega}\gamma\langle \nabla_Mu,\nabla_M\eta_d\rangle_Mu\eta_d\, d_M.
\end{multline*}
We obtain that
\begin{multline*}
\int_{\Omega_M^{d/2}\backslash\Omega_M^{3d/2}}\gamma\langle \nabla_Mu,\nabla_Mu\rangle_M\eta_d^2\, d_M=
\int_{\Omega}\gamma\langle \nabla_Mu,\nabla_Mu\rangle_M\eta_d^2\, d_M\\=
-2\int_{\Omega}\gamma\langle \nabla_Mu,\nabla_M\eta_d\rangle_Mu\eta_d\, d_M=
-2\int_{\Omega_M^{d/2}\backslash\Omega_M^{3d/2}}\gamma\langle \nabla_Mu,\nabla_M\eta_d\rangle_Mu\eta_d\, d_M
\\\leq \frac{2C}{d}\left(\int_{\Omega_M^{d/2}\backslash\Omega_M^{3d/2}}\gamma\langle \nabla_Mu,\nabla_Mu\rangle_M\eta_d^2\, d_M\right)^{1/2}\left(\int_{\Omega_M^{d/2}\backslash\Omega_M^{3d/2}}\gamma u^2\, d_M\right)^{1/2}
\end{multline*}
and we conclude that
\begin{equation}\label{aaa}
\int_{\Omega_M^{d/2}\backslash\Omega_M^{3d/2}}\gamma\langle \nabla_Mu,\nabla_Mu\rangle_M\eta_d^2\, d_M\leq
\frac{4C^2}{d^2}\int_{\Omega_M^{d/2}\backslash\Omega_M^{3d/2}}\gamma u^2\, d_M.
\end{equation}

Since
$$\int_{\Omega_M^{3d/4}\backslash \Omega_M^{5d/4}}\gamma\langle \nabla_Mu,\nabla_Mu\rangle_M\, d_M
\leq \int_{\Omega_M^{d/2}\backslash\Omega_M^{3d/2}}\gamma\langle \nabla_Mu,\nabla_Mu\rangle_M\eta_d^2\, d_M,
$$
by \eqref{aaa} and \eqref{bbb}, we infer that
\begin{equation}\label{ccc}
\int_{\Omega_M^{3d/4}\backslash \Omega_M^{5d/4}}\gamma\langle \nabla_Mu,\nabla_Mu\rangle_M\, d_M
\leq 
H(0)\frac{4C^2}{d}
e^{3C_3d/2}h(c_3d\Phi_1/2).
\end{equation}

Now we consider the function $u_d=u\eta_{d/4}$ and we easily prove that
\begin{multline*}
\int_{\Omega_M^d}\gamma\|\nabla_M u_d\|^2_M\, d_M\leq 2
\int_{\Omega_M^d}\gamma\left(\eta_{d/4}^2\|\nabla_M u\|^2_M+u^2\|\nabla_M\eta_{d/4}\|^2_M\right)\, d_M
\\\leq 
2\left(\int_{\Omega_M^{3d/4}\backslash \Omega_M^{5d/4}}\gamma\langle \nabla_Mu,\nabla_Mu\rangle_M\, d_M
+\frac{16C^2}{d^2}\int_{\Omega_M^{3d/4}\backslash \Omega_M^{5d/4}}\gamma u^2\, d_M\right)\\\leq
H(0)\frac{40C^2}{d}
e^{3C_3d/2}h(c_3d\Phi_1/2).
\end{multline*}

We have that $w_d=u-u_d$ solves, in a weak sense,
$$\left\{\begin{array}{ll}
\mathrm{div}_M(\gamma\nabla_M w_d)=-\mathrm{div}_M(\gamma\nabla_M u_d)&\text{in }\Omega_M^d\\
w_d=0&\text{on }\partial\Omega_M^d,
\end{array}
\right.
$$
from which we deduce that
$$\left(\int_{\Omega_M^d}\gamma\|\nabla_M w_d\|^2_M\, d_M\right)^{1/2}\leq\left(\int_{\Omega_M^d}\gamma\|\nabla_M u_d\|^2_M\, d_M\right)^{1/2},$$
hence
$$D(d)\leq H(0)\frac{160C^2}{d}
e^{3C_3d/2}h(c_3d\Phi_1/2).
$$
and the proof of \eqref{decayestimateter} is concluded by taking $C_4=160C^2$.\cvd

\smallskip

The rest of the section is devoted to the proofs of Theorems~\ref{mainthm} and \ref{mainthmbis}.

We also need the following notation. We note that, since $u\in H^2_{\mathrm{loc}}(\Omega)$, for any $d$ with $0<d<d_0$, $\nabla u$ is well-defined, in the trace sense, on 
$\partial\Omega^d_M$ and that $\nabla u\in L^2(\partial\Omega^d_M,\mathbb{R}^N)$.

For any $d$ with $0<d<d_0$, and almost any $x\in\partial\Omega^d_M$, with respect to the $(N-1)$-dimensional Hausdorff measure, we call $u_{\nu_M}(x)$ the (exterior) normal derivative of $u$ at $x$ with respect to $\Omega^d_M$ in the Riemannian setting which is given by
$$u_{\nu_M}(x)=\langle(\nabla_Mu(x))^T,\nu_M(x)\rangle_M=\nabla u(x) \nu_M(x).$$
We note that, analogously,
 $u_{\nu_M}\in L^2(\partial\Omega^d_M)$ is well-defined, again in the trace sense.
Moreover, using the equation and the divergence theorem, we have, for any $d$ with $0<d<d_0$,
$$D(d)=\int_{\partial\Omega^d_M}\gamma(x)u(x) u_{\nu_M}(x)\, d\sigma_M(x).$$
Finally we call, for any $d$ with $0<d<d_0$,
$$T(d)=\int_{\partial\Omega^d_M}\gamma(x)u^2_{\nu_M}(x)\, d\sigma_M(x)\quad\text{and}\quad F(d)=\frac{T(d)}{D(d)}.$$
We note that,
by a simple application of the Cauchy-Schwarz inequality, we have
$$F(d)\geq N(d)\qquad\text{for any }0<d<d_0.$$

Following essentially the arguments developed in \cite{Gar-Lin1}, we compute the derivatives, with respect to $d$, of $D$ and $H$.

By coarea formula and the properties of $\varphi_M$, we infer that $D$ is absolutely continuous on every compact subinterval contained in $[0,d_0)$ and that, for almost every $d\in (0, d_0)$,
$$D'(d)=-\int_{\partial\Omega^d_M}\gamma(x)\|\nabla_M u(x)\|_M^2\, d\sigma_M(x).$$
Then we use the following lemma, which is a suitable version of the Rellich identity.

\begin{lem}\label{Rellichid}
Let $u\in H^2_{\mathrm{loc}}(\Omega)$, $\gamma\in C^{0,1}(\overline{\Omega})$ and $v=(v^1,\ldots,v^N)\in C^{0,1}(\overline{\Omega},\mathbb{R}^N)$. Then
\begin{multline}\label{rellich}
\mathrm{div}_M(\gamma\langle\nabla_M u,\nabla_M u\rangle_M v)+2\mathrm{div}_M(\gamma \nabla_M u)\langle\nabla_M u,v\rangle_M
\\
=\gamma\langle\nabla_M u,\nabla_M u\rangle_M\mathrm{div}_M(v)+(\nabla\gamma v)\langle\nabla_M u,\nabla_M u\rangle_M+
\gamma(\nabla g^{l,j}v)u_lu_j\\+2\mathrm{div}_M(\gamma\langle\nabla_M u,v\rangle_M \nabla_M u)-2\gamma g^{l,j}u_lu_kv^k_j.
\end{multline}

\end{lem}

\proof{.} 
It follows by straightforward computations. In fact, with the summation convention and using subscripts for partial derivatives,
\begin{multline*}
\mathrm{div}_M(\gamma\langle\nabla_M u,\nabla_M u\rangle_M v)=
\gamma\langle\nabla_M u,\nabla_M u\rangle_M\mathrm{div}_M(v)+(\gamma g^{l,j}u_lu_j)_kv^k\\
=\gamma\langle\nabla_M u,\nabla_M u\rangle_M\mathrm{div}_M(v)+(\nabla\gamma v)(g^{l,j}u_lu_j)+
\gamma(\nabla g^{l,j}v)u_lu_j+\gamma g^{l,j}(u_{lk}u_j+u_lu_{jk})v^k.
\end{multline*}
On the other hand,
\begin{multline*}
2\mathrm{div}_M(\gamma\langle\nabla_M u,v\rangle_M \nabla_M u)=
2\mathrm{div}_M(\gamma \nabla_M u)\langle\nabla_M u,v\rangle_M+2\gamma \langle \nabla_M u,\nabla(\nabla uv)\rangle\\
=2\mathrm{div}_M(\gamma \nabla_M u)\langle\nabla_M u,v\rangle_M+2\gamma g^{l,j}u_l(u_kv^k)_j\\
=2\mathrm{div}_M(\gamma \nabla_M u)\langle\nabla_M u,v\rangle_M+2\gamma g^{l,j}u_lu_kv^k_j
+2\gamma g^{l,j}u_lu_{kj}v^k.
\end{multline*}
Finally,
$$\gamma g^{l,j}u_lu_{kj}v^k=\gamma g^{l,j}u_lu_{jk}v^k=\gamma g^{l,j}u_{lk}u_jv^k,$$
which follows by symmetry of the Hessian matrix and by observing that
$\gamma g^{l,j}u_{lk}u_jv^k=\gamma g^{j,l}u_{jk}u_lv^k=\gamma g^{l,j}u_lu_{jk}v^k$
since again by symmetry $g^{l,j}=g^{j,l}$.

Putting the three previous equality together, the lemma is proved.\cvd

\smallskip

We construct a Lipschitz function $v$ on $\overline{\Omega}$, with values in $\mathbb{R}^N$, coinciding with $\nu_M$ in $U_M^{d_0}$. By Remark~\ref{regUd} and Proposition~\ref{extensionprop},
or with a much simpler argument if \eqref{crucialequality2} holds, we can construct $v$ in such a way that 
$\|v\|_{C^{0,1}(\overline{\Omega},\mathbb{R}^N)}$ is bounded by a constant depending on the a priori data only.

Then we apply the Rellich identity \eqref{rellich} in $\Omega^d_M$, with $0<d<d_0$, to our solution $u$ and such a function $v$. Namely, by the divergence theorem,
\begin{multline*}
\int_{\partial\Omega^d_M}\gamma\langle\nabla_M u,\nabla_M u\rangle_M\, d\sigma_M=\int_{\partial\Omega^d_M}\gamma\langle\nabla_M u,\nabla_M u\rangle_M\langle v,\nu_M\rangle_M\, d\sigma_M\\= 
\int_{\Omega^d_M}(\mathrm{div}_M(\gamma\langle\nabla_M u,\nabla_M u\rangle_M v)+2\mathrm{div}_M(\gamma \nabla_M u)\langle\nabla_M u,v\rangle_M)\, d_M\\=
2\int_{\Omega^d_M}\mathrm{div}_M(\gamma\langle\nabla_M u,v\rangle_M \nabla_M u)\, d_M+A_0(d)=
2\int_{\partial\Omega^d_M}\gamma u^2_{\nu_M}\, d\sigma_M+A_0(d)
\end{multline*}
where
$$
A_0(d)=\int_{\Omega^d_M}\left(\langle\nabla_M u,\nabla_M u\rangle_M(\gamma\mathrm{div}_M(v)+\nabla\gamma v)+\gamma(\nabla g^{l,j}v)u_lu_j-2\gamma g^{l,j}u_lu_kv^k_j
\right)\, d_M.
$$
In other words, for almost every $d\in (0, d_0)$,
$$D'(d)=-2T(d)-A_0(d).$$
Finally, it is not difficult to show that there exists a positive constant $C$, depending on the a priori data only, such that for any $d$
with $0<d<d_0$ we have
$$|A_0(d)|\leq C D(d),$$
consequently, for almost every $d\in (0, d_0)$,
\begin{equation}\label{Dprime}
2F(d)-C\leq-\frac{D'(d)}{D(d)}\leq 2F(d)+C.
\end{equation}

Now we turn to the computation of $H'$. We wish to prove a similar estimate, namely that there exists a positive constant $\tilde{C}$, depending on the a priori data only, such that for almost every $d\in (0, d_0)$,
\begin{equation}\label{Hprime}
2N(d)-\tilde{C}\leq-\frac{H'(d)}{H(d)}\leq 2N(d)+\tilde{C}.
\end{equation}
Such a result directly follows by applying Proposition~\ref{normalderivativelemma}
to $w=\gamma u^2$.  In fact, we obtain that $H$ is absolutely continuous on every compact subinterval contained in $[0,d_0)$ and that,
for almost any $d$, $0<d<d_0$,
$$-H'(d)=2D(d)+\int_{\partial\Omega^d_M}(\nabla \gamma(x)\nu_M(x))u^2(x)\, d\sigma_M(x)+A(d)=2D(d)+A_1(d).$$
Again, it is not difficult to show that there exists a positive constant $\tilde{C}$, depending on the a priori data only, such that for any $d$
with $0<d<d_0$ we have
$$|A_1(d)|\leq \tilde{C} H(d),$$
consequently, for almost every $d\in (0, d_0)$, \eqref{Hprime} holds.

We are now in the position to conclude the proof of Theorem~\ref{mainthm}.

\smallskip

\proof{ of Theorem~\textnormal{\ref{mainthm}}.}
We use an ordinary differential equation argument, exploiting \eqref{Dprime} and \eqref{Hprime}.
With $C$ as in \eqref{Dprime}, let us define, for $0\leq d<d_0$,
$$\tilde{D}(d)=e^{-Cd}D(d).$$
Then
$$-\frac{\tilde{D}'(d)}{\tilde{D}(d)}\geq 2F(d)=2(F(d)-N(d))+2N(d).$$
Therefore, for any $0<d<d_0$,
$$\log(\tilde{D}(0))-\log(\tilde{D}(d))\geq 2\int_0^d(F(t)-N(t))\, dt+2\int_0^d N(t)\, dt.$$
We call
$$G_0(d)=e^{-2\int_0^dF}\quad\text{and}\quad G(d)=e^{-2\int_0^d(F-N)}\quad\text{and}\quad G_1(d)=e^{-2\int_0^dN}$$
and we obtain that
$$\tilde{D}(d)\leq G_0(d)\tilde{D}(0)=G(d)G_1(d)\tilde{D}(0),$$
therefore
\begin{equation}\label{firsteq}
D(d)\leq e^{Cd}G_0(d)D(0)=e^{Cd}G(d)G_1(d)D(0).
\end{equation}

Since
$$\frac{N'(d)}{N(d)}=\frac{D'(d)}{D(d)}-\frac{H'(d)}{H(d)},$$
we infer that
\begin{equation}\label{logNder}
2(F(d)-N(d))-\hat{C}\leq-\frac{N'(d)}{N(d)}\leq 2(F(d)-N(d))+\hat{C}
\end{equation}
where $\hat{C}=C+\tilde{C}$. Let us define, for $0\leq d<d_0$,
$$\tilde{N}(d)=e^{-\hat{C}d}N(d).$$
Then,
by \eqref{logNder}, we conclude that
$$-\frac{\tilde{N}'(d)}{\tilde{N}(d)}\geq 2(F(d)-N(d))\geq 0.$$
In other words, $\tilde{N}$ is decreasing. We note that this is the crucial point in the argument of \cite{Gar-Lin1}.
However, in our case, such a property is not enough, since, in order to estimate $G_1$, we need to control how fast $N$ can decrease. Still by \eqref{logNder}, we infer that
\begin{equation}\label{secondeq}
N(d)\geq e^{-\hat{C}d}G(d)N(0).
\end{equation}

We note that $G(0)=1$ and, since $F-N\geq 0$, $G$ is positive and decreasing with respect to $d$. We estimate $G_1(d)$ by using \eqref{secondeq} and the fact that $G(s)\geq G(d)$ for any $0<s<d$, obtaining that
$$G_1(d)\leq e^{-2N(0)G(d)\int_0^de^{-\hat{C}s}\, ds}=
e^{-b(d)N(0)G(d)},$$
where
$$b(d)=\frac{2}{\hat{C}}\left(1-e^{-\hat{C}d}\right).$$

We consider the auxiliary function $g(x)=xe^{-\alpha x}$, $x\in[0,1]$, with $\alpha>0$, and note that
$$\max_{x\in[0,1]}g(x)=h(\alpha),$$
thus
we conclude that
\begin{equation}\label{finaleq}
D(d)\leq e^{Cd}D(0)h(b(d)N(0)).
\end{equation}

Since $\lambda_1c_1\Phi\leq N(0)$ and $2e^{-\hat{C}d_0}d\leq b(d)$, the proof of \eqref{decayestimate} is concluded
by setting $C_2=C$ and $c_2=2\lambda_1c_1e^{-\hat{C}d_0}$.\cvd

\smallskip

We note that, without any control on $F-N$, besides the fact that it is positive, using this technique it is in practice impossible to improve the estimate of
Theorem~\ref{mainthm}.

%
%

We now turn to the proof of Theorem~\ref{mainthmbis}. We need the following notation.
For any $d\in [0,d_0)$ we define
$$E(d)=\int_{\Omega^d_M}\gamma u^2\, d_M.$$
We note that $E$ is a strictly positive function which is absolutely continuous on any compact subinterval of $[0,d_0)$. Moreover, for almost any $d\in(0,d_0)$, we have
$$E'(d)=-\int_{\partial\Omega^d_M}\gamma u^2\, d\sigma_M=-H(d).$$

We construct a Lipschitz function $v_1\in C^{0,1}(\overline{\Omega},\mathbb{R}^N)$ coinciding with $\nu_M$ in $U_M^{d_0/2}$ and such that
$\|v_1\|_{C^{0,1}(\overline{\Omega},\mathbb{R}^N)}$ is bounded by a constant depending on the a priori data only and
$$\|v_1\|_M\leq 1\quad\text{in }\overline{\Omega}.$$
Such a construction is fairly easy. We consider a $C^{\infty}$ function $\chi:\mathbb{R}\to\mathbb{R}$ such that
$\chi$ is increasing, $\chi=0$ on $(-\infty,3d_0/5]$ and $\chi=1$ on $[4d_0/5,+\infty)$. We can define $v_1$ with the desired properties as follows
$$v_1(x)=\chi(\varphi_M(x))\frac{e_1}{\sqrt{\langle G(x)e_1,e_1\rangle}}+(1-\chi(\varphi_M(x)))\nu_M(x)\quad\text{for any }x\in\overline{\Omega}.$$

Then we have, for any $d$, $0\leq d\leq d_0/2$,
\begin{multline}\label{Sdef}
H(d)=
\int_{\partial\Omega^d_M}\gamma u^2\, d\sigma_M=
\int_{\Omega^d_M}\mathrm{div}_M(\gamma u^2v_1)\, d_M\\=
2\int_{\Omega^d_M}\gamma u\langle(\nabla_M u)^T,v_1\rangle_M\, d_M+\int_{\Omega^d_M}\mathrm{div}_M(\gamma v_1)u^2\, d_M=2S(d)+A_2(d).
\end{multline}
It is not difficult to show that, for some constant $\tilde{C}_1$ depending on the a priori data only, we have
\begin{equation}\label{Sdefbound}
|A_2(d)|\leq\tilde{C}_1E(d).
\end{equation}

We now call, for any $d\in [0,d_0)$,
$$K(d)=\frac{H(d)}{E(d)}\quad\text{and}\quad K_1(d)=\frac{H(d)}{\sqrt{E(d)}}.$$
For almost any $d$ with $0<d<d_0$, we have
\begin{equation}\label{Eprime}
-\frac{E'(d)}{E(d)}=K(d).
\end{equation}
Since
$$\frac{K'(d)}{K(d)}=\frac{H'(d)}{H(d)}-\frac{E'(d)}{E(d)},$$
by \eqref{Hprime} and \eqref{Eprime}
we infer that, for almost every $d\in (0,d_0)$,
\begin{equation}\label{logKder}
2N(d)-K(d)-\tilde{C}\leq-\frac{K'(d)}{K(d)}\leq 2N(d)-K(d)+\tilde{C}.
\end{equation}
Analogously, since
$$\frac{K_1'(d)}{K_1(d)}=\frac{H'(d)}{H(d)}-\frac{1}{2}\frac{E'(d)}{E(d)},$$
we obtain that, for almost every $d\in (0,d_0)$,
\begin{equation}\label{logK1der}
2N(d)-\frac{K(d)}{2}-\tilde{C}\leq-\frac{K_1'(d)}{K_1(d)}\leq 2N(d)-\frac{K(d)}{2}+\tilde{C}.
\end{equation}

We are now in the position to conclude the proof of Theorem~\ref{mainthmbis}.

\smallskip

\proof{ of Theorem~\textnormal{\ref{mainthmbis}}.}
In the sequel, we adopt the following normalisation, that is, we assume that
\begin{equation}\label{Enorm}
E(0)=1.
\end{equation}
It is immediate to show, with this assumption, that for any $d\in [0,d_0)$ we have
$E(d)\leq 1$ and, consequently,
\begin{equation}\label{KvsK1}
K_1(d)\leq K(d).
\end{equation}

We now apply a similar technique we used before to estimate $D$ to the function $H$.
Namely,
with $\tilde{C}$ as in \eqref{Hprime}, let us define, for any $d\in [0,d_0)$,
$$\tilde{H}(d)=e^{-\tilde{C}d}H(d).$$
Then, by \eqref{KvsK1},
$$-\frac{\tilde{H}'(d)}{\tilde{H}(d)}\geq 2N(d)\geq \left(2N(d)-\frac{K(d)}{2}\right)+\frac{K_1(d)}{2}.$$
The main difference with respect to the argument for $D$ is that, whereas it is immediate to show that $F\geq N$,
it is not that evident that $2N\geq K/2$. The other difference with respect to the previous argument is that we need to use $K_1$ instead of $K$ itself. However, for any $d\in[0,d_0/2)$, using \eqref{Sdef} and calling $\tilde{A}_2(d)=A_2(d)/2$,
\begin{multline*}
2N(d)-\frac{K(d)}{2}=2N(d)-\frac{S(d)}{E(d)}-\frac{\tilde{A}_2(d)}{E(d)}=
\frac{D(d)}{S(d)+\tilde{A}_2(d)}-\frac{S(d)}{E(d)}-\frac{\tilde{A}_2(d)}{E(d)}\\
=2\frac{D(d)E(d)-S^2(d)-S(d)\tilde{A}_2(d)}{H(d)E(d)}-\frac{\tilde{A}_2(d)}{E(d)}=
2\frac{D(d)E(d)-S^2(d)}{H(d)E(d)}-\frac{S(d)A_2(d)}{H(d)E(d)}-\frac{\tilde{A}_2(d)}{E(d)}.
\end{multline*}
It is easy to see that the first term
$$M(d)=2\frac{D(d)E(d)-S^2(d)}{H(d)E(d)}$$
is positive, therefore, using \eqref{Sdef} and \eqref{Sdefbound}, we obtain that
\begin{multline}\label{aneq}
2N(d)-\frac{K(d)}{2}= M(d)-\frac{A_2(d)}{E(d)}\left[1-\frac{1}{2}\frac{A_2(d)}{H(d)}\right]\\
= M(d)+\frac{A_2^2(d)}{2H(d)E(d)} -\frac{A_2(d)}{E(d)}
=M_1(d)-\frac{A_2(d)}{E(d)}
\end{multline}
where
\begin{equation}\label{M1defin}
M_1(d)=M(d)+\frac{A_2^2(d)}{2H(d)E(d)}\geq 0\quad\text{for any }d\in (0,d_0/2).
\end{equation}
We note that, by \eqref{Sdefbound}, for any $d$ with $0<d<d_0/2$
\begin{equation}\label{rightone}
M_1(d)-\tilde{C}_1\leq 2N(d)-\frac{K(d)}{2}\leq M_1(d)+\tilde{C}_1,
\end{equation}
consequently, by \eqref{logK1der} and calling $\tilde{C}_2=\tilde{C}+\tilde{C}_1$, we have,
for almost any $d\in (0,d_0/2)$,

\begin{equation}\label{rightonebis}
M_1(d)-\tilde{C}_2\leq -\frac{K_1'(d)}{K_1(d)}\leq M_1(d)+\tilde{C}_2.
\end{equation}

For any $d$ with $0<d<d_0/2$, we have
$$\log(\tilde{H}(0))-\log(\tilde{H}(d))\geq \int_0^dM_1(t)\, dt+\frac{1}{2}\int_0^d K_1(t)\, dt-\tilde{C}_1d.$$
We call
$$J_0(d)=e^{-2\int_0^dN}\quad\text{and}\quad J(d)=e^{-\int_0^dM_1}\quad\text{and}\quad J_1(d)=e^{-\int_0^d(K_1/2)}$$
and we obtain that
$$\tilde{H}(d)\leq J_0(d)\tilde{H}(0)\leq e^{\tilde{C}_1d}J(d)J_1(d)\tilde{H}(0),$$
therefore,
\begin{equation}\label{firsteqH}
H(d)\leq e^{\tilde{C}d}J_0(d)H(0)\leq e^{\tilde{C}_2d}J(d)J_1(d)H(0).
\end{equation}

For any $d$ with $0\leq d<d_0/2$, by \eqref{rightonebis}
and since $K_1(0)=K(0)$,
\begin{equation}\label{secondeqbis}
K_1(d)\geq e^{-\tilde{C}_2d}J(d)K(0).
\end{equation}

We note that $J(0)=1$ and, since $M_1\geq 0$, $J$ is positive and decreasing with respect to $d$. We estimate $J_1(d)$ by using \eqref{secondeqbis} and the fact that $J(s)\geq J(d)$ for any $0<s<d$, obtaining that
$$J_1(d)\leq e^{-(K(0)/2)J(d)\int_0^de^{-\tilde{C}_2s}\, ds}=
e^{-\tilde{b}(d)K(0)J(d)},$$
where
$$\tilde{b}(d)=\frac{1}{2\tilde{C}_2}\left(1-e^{-\tilde{C}_2d}\right).$$

Arguing as in the proof of Theorem~\ref{mainthm}, we conclude
that, setting $C_3=\tilde{C}_2$, for any $d$ with $0<d<d_0/2$ we have
\begin{equation}
H(d)\leq e^{C_3d}H(0)h(\tilde{b}(d)K(0)).
\end{equation}

In order to conclude the proof it is enough to show that, for some positive constant $c_3$ depending on the a priori data only,
we have for any $d$ with $0<d<d_0/2$
$$\tilde{b}(d)K(0)\geq c_3d\Phi_1.$$
This is an immediate consequence of Proposition~\ref{-1/2-2bound}.\cvd

\end{document}